\documentclass[amscd,amssymb]{amsart}

\usepackage[all]{xy}

\newtheorem{thm}{Theorem}[section]
\newtheorem{lem}[thm]{Lemma}
\newtheorem{cor}[thm]{Corollary}
\newtheorem{prop}[thm]{Proposition}

\theoremstyle{definition}

\newtheorem{defn}[thm]{Definition}

\newtheorem{ex}[thm]{Example}

\newtheorem{eqn}[thm]{}

\theoremstyle{remark}

\newtheorem{rem}[thm]{Remark}

\numberwithin{equation}{section}

\newcommand{\thmref}[1]{Theorem~\ref{#1}}
\newcommand{\corref}[1]{Corollary~\ref{#1}}
\newcommand{\secref}[1]{\S\ref{#1}}

\newcommand{\propref}[1]{Proposition~\ref{#1}}
\newcommand{\lemref}[1]{Lemma~\ref{#1}}

\newcommand{\exref}[1]{Example~\ref{#1}}
\newcommand{\remref}[1]{Remark~\ref{#1}}

\newcommand{\Hom}{\operatorname{Hom}}
\newcommand{\Ext}{\operatorname{Ext}}

\newcommand{\Ind}{\operatorname{Ind}}
\newcommand{\Res}{\operatorname{Res}}
\newcommand{\Rep}{\operatorname{Rep}}
\newcommand{\Epi}{\operatorname{Epi}}
\newcommand{\rk}{\operatorname{rk}}
\newcommand{\mrk}{\operatorname{mrk}}
\newcommand{\ob}{\operatorname{ob}}
\newcommand{\depth}{\operatorname{depth}}
\newcommand{\Ess}{\operatorname{Ess}}

\newcommand{\A}{{\mathcal  A}}

\newcommand{\U}{{\mathcal  U}}
\newcommand{\K}{{\mathcal  K}}
\newcommand{\HH}{{\mathcal  H}}

\newcommand{\Z}{{\mathbb  Z}}

\newcommand{\F}{{\mathbb  F}}

\newcommand{\ra}{\rightarrow}
\newcommand{\xra}{\xrightarrow}
\newcommand{\la}{\leftarrow}
\newcommand{\xla}{\xleftarrow}

\newcommand{\lra}{\longrightarrow}

\begin{document}

\title[group automorphisms]{The nilpotent filtration and the action of automorphisms on the cohomology of finite $p$--groups}

 \author[Kuhn]{Nicholas J.~Kuhn}
 \address{Department of Mathematics \\ University of Virginia \\ Charlottesville, VA 22904}
 \email{njk4x@virginia.edu}
\thanks{This research was partially supported by a grant from the National Science Foundation}

 \date{May 3, 2006.}


 \begin{abstract}
We study $H^*(P)$, the mod $p$ cohomology of a finite $p$--group $P$, viewed as an $\F_p[Out(P)]$--module.  In particular, we study the conjecture, first considered by Martino and Priddy, that, if $e_S \in \F_p[Out(P)]$ is a primitive idempotent associated to an irreducible $\F_p[Out(P)]$--module $S$, then the Krull dimension of $e_SH^*(P)$ equals the rank of $P$. The rank is an upper bound by Quillen's work, and the conjecture can be viewed as the statement that every irreducible $\F_p[Out(P)]$--module occurs as a composition factor in $H^*(P)$ with similar frequency.

In summary, our results are as follows.  A strong form of the conjecture is true when $p$ is odd.  The situation is much more complex when $p=2$, but is reduced to a question about 2--central groups (groups in which all elements of order 2 are central), making it easy to verify the conjecture for many finite 2--groups, including all groups of order 128, and all groups that can be written as the product of groups of order 64 or less.

The the odd prime theorem can be deduced using the approach to $\U$, the category of unstable modules over the Steenrod algebra, initiated by H.-W. Henn, J. Lannes, and L. Schwartz in \cite{hls1}.  The reductions when $p=2$ make heavy use of the nilpotent filtration of $\U$ introduced in \cite{s1}, as applied to group cohomology in \cite{hls2}.  Also featured are unstable algebras of cohomology primitives associated to central group extensions.

\end{abstract}

\maketitle

\section{Introduction}

Fix a prime $p$, and let $H^*(P)$ denote the mod $p$ cohomology ring of a finite $p$--group $P$. Interpreted topologically, $H^*(P)=H^*(BP;\F_p)$, where $BP$ is the classifying space of $P$, and interpreted algebraically, $H^*(P) = \Ext_{\F_p[P]}^*(\F_p,\F_p)$.

The automorphism group of $P$, $Aut(P)$, acts on $P$, and thus via ring homomorphisms on $H^*(P)$.  As the inner automorphism group, $Inn(P)$, acts trivially on cohomology, $H^*(P)$ becomes a graded $\F_p[Out(P)]$--module, where $Out(P)=Aut(P)/Inn(P)$ is the outer automorphism group.

The ring $H^*(P)$ is known to be Noetherian, and D. Quillen \cite{quillen} computed its Krull dimension:
$ \dim H^*(P) = \rk(P)$.  Here $\rk(P)$ denotes the maximal rank of an elementary abelian $p$--subgroup of $P$.

Now let $S$ be an irreducible $\F_p[Out(P)]$--module and $e_S \in \F_p[Out(P)]$ an associated primitive idempotent.  Then $e_SH^*(P)$ is a finitely generated module over the Noetherian ring $H^*(P)^{Out(P)}$ (see \secref{finite subsection}), and thus has a Krull dimension with evident upper bound $\rk(P)$.  In the early 1990's, J. Martino and S. Priddy \cite{mp} asked whether the following conjecture might be true. \\

\noindent{\bf Conjecture A} \ $\dim e_SH^*(P) = \rk(P)$ for all pairs $(P,S)$. \\

As $\dim e_SH^*(P)$ determines the growth of the Poincar\'e series of $e_SH^*(P)$, and this Poincar\'e series gives a count of the occurrences of $S$ as a composition factor in $H^*(P)$, the conjecture is roughly the statement that every irreducible $\F_p[Out(P)]$--module occurs in $H^*(P)$ with similar frequency.

The topological interpretation of $H^*(P)$ makes it evident that $H^*(P)$ is an object in $\K$ and $\U$, the categories of unstable algebras and modules over the mod $p$ Steenrod algebra $\A$, and in this paper, we study Conjecture A using modern `$\U$--technology'.

In summary, our results are as follows.  A strong form of Conjecture A is true when $p$ is odd.  The situation is more complex when $p=2$.  We reduce the conjecture to a different conjecture about 2--central groups, making it easy to verify Conjecture A for many finite 2--groups, including all with order dividing 128, and all groups that can be written as the product of groups of order 64 or less.

The the odd prime theorem can be deduced using the functor category approach to $\U$ initiated by H.-W. Henn, J. Lannes, and L. Schwartz in \cite{hls1}.  The reductions when $p=2$ make heavy use of the nilpotent filtration of $\U$ introduced in \cite{s1}, as applied to group cohomology in \cite{hls2}.

Featured here, and in a companion paper \cite{k4}, are algebras of cohomology primitives associated to central group extensions, and  $p$--central groups: groups in which every element of order $p$ is central.

\section{Main results}

\subsection{Notation} Before describing our results in more detail, we introduce some notation. Throughout, $V < P$ will denote an elementary abelian $p$--subgroup of $P$ and $C = C(P) < P$ will be the maximal central elementary abelian $p$--subgroup.  We let $c(P)$ be the rank of $C$, $\depth_{H^*(P)}H^*(P)$ be the depth of $H^*(P)$ with respect to the ideal of positive degree elements, and $\mrk(P)$ be the minimal rank of a maximal $V < P$.
These are related by
$$c(P) \leq \depth_{H^*(P)}H^*(P) \leq \mrk(P) \leq \rk(P),$$ where the first inequality is due to J.Duflot \cite{duflot}, and the second follows from Quillen's work and standard commutative algebra (see \secref{depth subsection}).

Given $V < P$, we let $Aut_0(P,V) \lhd Aut(P,V) < Aut(P)$ be defined by
$$ Aut(P,V) = \{ \alpha \in Aut(P) \ | \ \alpha(V) = V \} \text{ and} $$
$$  Aut_0(P,V) = \{ \alpha \in Aut(P) \ | \ \alpha(v) = v \text{ for all } v \in V \},$$
Then let $Out_0(P,V)$ and $Out(P,V)$ be the projection of these groups in $Out(P)$.
Note that $Out(P,V) = Aut(P,V)/N_P(V)$ and $Out_0(P,V) = Aut_0(P,V)/C_P(V)$, so that there is an extension of groups
$$ W_P(V) \ra  Aut(P,V)/Aut_0(P,V) \ra Out(P,V)/Out_0(P,V),$$
where $W_P(V) = N_P(V)/C_P(V)$.

Finally, if $S$ is an irreducible $\F_p[G]$--module, and $M$ is an arbritrary $\F_p[G]$--module, we say that {\em $S$ occurs in $M$} if it is a composition factor. \\

\subsection{Quillen's approximation}

Let $\A(P)$ denote Quillen's category, having as objects the elementary abelian $p$--subgroups $V < P$, and as morphisms the homomorphisms generated by inclusions and conjugation by elements in $P$.  Let $\HH^*(P)$ denote $\displaystyle \HH^*(P) = \lim_{\A(P)} H^*(V)$.  It is not hard to see the same inverse limit is attained by using the smaller category $\A^C(P)$, where $\A^C(P) \subset \A(P)$ is the full subcategory consisting of elementary abelian $V < P$ containing the maximal central elementary abelian $C$.

The restriction maps associated to the inclusions $V < P$ assemble to define a natural map of unstable $\A$--algebras
$$ H^*(P) \ra \HH^*(P),$$
which Quillen shows is an $F$--isomorphism.  As a special case of \propref{Rd prop} below, we have
\begin{equation*}  \dim e_SH^*(P) \geq \dim e_S\HH^*(P),
\end{equation*}
for all irreducible $\F_p[Out(P)]$--modules $S$.

\begin{thm}  \label{R_0 theorem}  Let $S$ be an irreducible $\F_p[Out(P)]$--module.  Then $e_S\HH^*(P) \neq 0$ if and only if $S$ occurs in $\F_p[Out(P)/Out_0(P,V)]$ for some maximal $V < P$.  In that case,
$\dim e_S\HH^*(P) = \max \{  \rk(V)  \ | \ S \text{ occurs in } \F_p[Out(P)/Out_0(P,V)] \}.$
\end{thm}

This implies an observation noted earlier by Martino and Priddy \cite[Prop.4.2]{mp}.

\begin{cor} \label{martino priddy cor} If $e_S\HH^*(P) \neq 0$, then $\dim e_S\HH^*(P) \geq \mrk(P)$.
\end{cor}

Given $H < G$, a standard argument shows that if $H$ is a $p$--group, then every irreducible $\F_p[G]$--module occurs as a submodule of $\F_p[G/H]$.  Thus the theorem implies that if $Out_0(P,V)$ is a $p$--group, then $\dim e_S\HH^*(P) \geq \rk(V)$ for all irreducibles $S$.  The converse will be true if $P$ has a unique maximal elementary abelian subgroup, for example, if $P$ is $p$--central: see \corref{little cor}.

One can now make effective use of group theory results about $p^{\prime}$--automorphisms of $p$--groups, particularly the Thompson $A \times B$ lemma \cite[Thm.5.3.4]{gor}.

\begin{thm}  \label{Out(P,V) thm} Let $P$ a finite $p$--group.  If $V<P$ is an elementary abelian subgroup, then $Out_0(P,V)$ will be a $p$--group if $Out_0(C_P(V),V)$ is a $p$--group.  If $V$ is also maximal, and $p$ is odd, this will always be the case.  Thus if $p$ is odd, $\dim e_S\HH^*(P) = \rk(P)$ for all irreducible $\F_p[Out(P)]$--modules $S$.
\end{thm}

For $2$--groups, one does not need to look far to see that story is quite different.  The quaternionic group of order 8, $Q_8$, is $2$--central with center $C$ of rank 1.  As $Aut(Q_8)$ clearly must act trivially on the center, $Out_0(Q_8,C) = Out(Q_8) \simeq \Sigma_3$, which is not a 2--group.

This example generalizes: it is the case $t=1$ in the following.

\begin{ex} \label{2C examples} Let $G_t$ be the $2$--Sylow subgroup of $SU_3(\F_{2^t})$.  This group is 2--central with center $C_t \simeq \F_{2^t}$, and $G_t/C_t \simeq \F_{2^{2t}}$.  The group $Out(G_t,C_t)$ contains a cyclic subgroup of order $2^t+1$.  Thus there is at least one irreducible $\F_2[Out(G_t)]$--module $S$ with $e_S\HH^*(G_t) = 0$.
\end{ex}

The next proposition, an easy application of properties of Lannes $T_V$--functor, shows that \thmref{Out(P,V) thm} can be applied to 2--groups with reduced cohomology, i.e., groups $P$ for which $H^*(P) \ra \HH^*(P)$ is monic.  This includes groups built up by iterated products and wreath products of $\Z/2$'s.

\begin{prop} \label{reduced prop} If $H^*(P)$ is reduced, and $V < P$ is a maximal elementary abelian subgroup, then $C_P(V) = V$, and thus $Out_0(C_P(V),V)$ is the trivial group.
\end{prop}

To say more about $\dim e_SH^*(P)$, we now begin to take a deeper look at $H^*(P)$ using the nilpotent filtration.

\subsection{A stratification of the problem}

In \secref{U section}, we study functors $\bar R_d: \U \ra \U$.  In brief, they are defined as follows.  An unstable module $M$ has a natural nilpotent filtration
$$ \dots \subset nil_2M \subset nil_1M \subset nil_0 M= M,$$
and $ nil_dM/nil_{d+1}M = \Sigma^d R_dM$,
where $R_dM$ is a reduced unstable module.  Then $\bar R_dM$ is defined as the nilclosure of $R_dM$.

Henn \cite{henn} has shown that the nilpotent filtration of a Noetherian unstable algebra is finite, so only finitely many of the the modules $\bar R_dH^*(P)$ will be nonzero.  Furthemore, $\bar R_0 H^*(P) = \HH^*(P)$, and each $\bar R_d H^*(P)$ is a finitely generated $\HH^*(P)$--module. The study of $\dim e_SH^*(P)$ stratifies as follows.

\begin{prop} \label{Rd prop} Let $S$ be an irreducible $\F[Out(P)]$--module.  Then
$$\dim e_SH^*(P) = \max_{d} \{ \dim e_S\bar R_dH^*(P) \}.$$
\end{prop}

This will be a special case of a more general statement about $\K$: see \propref{general Rd prop}.

Martino and Priddy observe \cite[Prop.4.1]{mp} that the Depth Conjecture, subsequently proved by D. Bourguiba and S. Zarati \cite{bour zarati}, implies that $\dim e_SH^*(P) \geq \depth_{H^*(P)}H^*(P)$.  We prove a stratified variant of this, and note that an analogue of Duflot's theorem holds.

\begin{prop} \label{depth prop} If $e_S\bar R_dH^*(P)\neq 0$, then
$$\dim e_S\bar R_dH^*(P) \geq \depth_{\HH^*(P)}\bar R_dH^*(P) \geq c(P).$$
\end{prop}

That $\dim e_S\bar R_dH^*(P) \geq c(P)$ is also a corollary of part (c) of \thmref{R_d theorem} below.

\subsection{Computing $\dim e_S\bar R_dH^*(P)$ }

In \cite{k4}, the work of \cite{hls2} will allow us to write down useful formulae for $\bar R_dH^*(P)$ analogous to
$$ \bar R_0H^*(P) = \lim_{V \in \A^C(P)}H^*(V).$$
Using these, we can generalize \thmref{R_0 theorem} to a statement about $\dim e_S\bar R_dH^*(P)$ for all $d$.

To state this, we need to define the primitives associated to a central extension.

The cohomology of an elementary abelian $p$--group is a Hopf algebra.  If $Q$ is a finite group, and $V < Q$ is a central elementary abelian $p$--subgroup, then multiplication $m: V \times Q \ra Q$ is a homomorphism, and the induced map in cohomology,
$  m^*: H^*(Q) \ra H^*(V \times Q) = H^*(V) \otimes H^*(Q)$
makes $H^*(Q)$ into an $H^*(V)$--comodule.  We then let $P_VH^*(Q)$ denote the primitives:
\begin{equation*}
\begin{split}
P_VH^*(Q) & = \{ x \in H^*(Q) \ | \ m^*(x) = 1 \otimes x \} \\
& = \text{Eq } \{ H^*(Q)
\begin{array}{c} m^* \\[-.08in] \longrightarrow \\[-.1in] \longrightarrow \\[-.1in] \pi^*
\end{array}
H^*(V \times Q) \},
\end{split}
\end{equation*}
where $\pi: V \times Q \ra Q$ is the projection.

As equalizers of algebra maps are algebras, $P_VH^*(Q)$ is again an unstable $\A$--algebra.  Note also that $P_VH^0(Q) = \F_p$.  It is not hard to check that $P_VH^1(Q) \simeq H^1(Q/V) = \Hom(Q/V, \F_p)$, via the inflation map.  The reader may find it illuminating to know that $P_VH^*(Q)$ is again Noetherian, and has Krull dimension equal to $\rk(Q) - \rk(V)$: see \cite{k3,k4}.

Our theorem about $\dim e_S\bar R_dH^*(P)$ now goes as follows.

\begin{thm} \label{R_d theorem}  Let $S$ be an irreducible $\F_p[Out(P)]$--module.  For each $d$, there exists a set $\sup_d(S) \subset \text{ob }\A^C(P)$, the {\em $d$--support of $S$}, with the following properties. \\

\noindent{\bf (a)} \ $\sup_d(S)$ is a union of $Aut(P)$--orbits in $\text{ob }\A^C(P)$. \\

\noindent{\bf (b)} \ Suppose that $V_1 < V_2 < P$, and $P_{V_1}H^d(C_P(V_1)) \ra P_{V_1}H^d(C_P(V_2))$ is monic.  Then $V_1 \in \sup_d(S)$ implies that $V_2 \in \sup_d(S)$. \\

\noindent{\bf (c)} \ $e_S\bar R_dH^*(P) \neq 0$ if and only if $sup_d(S) \neq \emptyset$, and, in this case, $$\dim e_S\bar R_dH^*(P) = \max \{\rk(V) \ | \ V \in {\sup}_d(S) \}.$$

\noindent{\bf (d)} \ Let $V \in \text{ob }\A^C(P)$ be maximal.  Then $V \in \sup_d(S)$ if and only if $S$ occurs in $\displaystyle \Ind_{Out(P,V)}^{Out(P)}([\F_p[Aut(P,V)/Aut_0(P,V)] \otimes P_VH^d(C_P(V))]^{W_P(V)})$.

\end{thm}

When $d=0$, we recover \thmref{R_0 theorem}.  In this case, the hypothesis in part (b) always holds.  Thus maximal elements in $sup_0(S)$ will always be maximal in $\text{ob }\A^C(P)$.  Thus part (d) implies that $V$ will be a maximal element in $\sup_0(S)$ if and only if $S$ occurs in
$$  \Ind_{Out(P,V)}^{Out(P)}(\F_p[Aut(P,V)/Aut_0(P,V)]^{W_P(V)})$$
which rewrites as
$$\F_p[Out(P,V)/Out_0(P,V)].$$
Thus part (c) implies the calculation of $\dim e_S\HH^*(P)$ in \thmref{R_0 theorem}.

\thmref{R_d theorem} now lets us recast Conjecture A as follows.

\begin{cor} \label{conj A equivalent} Conjecture A is true for a pair $(P,S)$ if and only if there exists a $V < P$ of maximal rank such that $S$ occurs in $$\displaystyle \Ind_{Out(P,V)}^{Out(P)}([\F_p[Aut(P,V)/Aut_0(P,V)] \otimes P_VH^*(C_P(V))]^{W_P(V)}).$$
\end{cor}

Starting from this, it is not hard to show (see \lemref{B implies A lemma}) that Conjecture A for a $2$--group $P$ would be implied by \\

\noindent{\bf Conjecture B} Let $P$ be a finite 2--group.  If $V < P$ is a maximal elementary abelian subgroup, then every irreducible $\F_2[Out_0(P,V)]$--module occurs in $P_VH^*(C_P(V))^{W_P(V)}$. \\

Note that the subgroups $C_P(V)$ that arise here are $2$--central.  Our general theory shows

\begin{prop} \label{pC prop} If $Q$ is $2$--central with maximal elementary abelian subgroup $C$, then Conjecture B is true for $Q$, i.e. every irreducible $\F_2[Out_0(Q,C)]$--module occurs as a composition factor in $P_CH^*(Q)$.
\end{prop}

We conjecture a subtle strengthening of this. \\

\noindent{\bf Conjecture C}  If $Q$ is $2$--central with maximal elementary abelian subgroup $C$, then every irreducible $\F_2[Out_0(Q,C)]$--module occurs as a submodule of $P_CH^*(Q)$.  \\

Again using the Thompson $A \times B$ lemma, we have

\begin{thm} \label{conj C implies conj B thm} Conjecture $B$ is true for the pair $(P,V)$ if the $2$--central group $C_P(V)$ satisfies Conjecture C.  Thus if $P$ is a 2--group such that $C_P(V)$ satisfies Conjecture C for some $V<P$ of maximal rank, then $\dim e_SH^*(P) = \rk(P)$ for all irreducible $\F_2[Out(P)]$--modules $S$.
\end{thm}

Though $Out_0(Q,C)$ need not be a 2--group if $Q$ is 2--central, its $2^{\prime}$ part tends to be very small, even in cases when $Out(Q)$ is quite complicated.  This makes it not hard to check the following theorem, using information about 2--groups of order dividing 64 available at the website \cite{carlson website}, or, in book form, \cite{carlson et al}.

\begin{thm} \label{2 central conj thm} Conjecture C is true for all 2--central groups that can be written as the product of groups of order dividing 64.
\end{thm}

A counterexample to Conjecture A would have to be a 2--group $P$ that is not 2--central, and having the property that the proper subgroups $C_P(V)$, with $V$ of maximal rank, are all counterexamples to Conjecture C.  Thus the last theorem implies

\begin{cor} Conjecture A is true for all 2--groups of order dividing 128, and all 2--groups that can be written as the product of groups of order 64 or less.
\end{cor}

We end this section with an example\footnote{Checking Conjecture A for this group led us to the formulation and proof of \thmref{R_d theorem}.} that illustrates the sorts of patterns that the numbers $\dim e_S\bar R_dH^*(P)$ can take, when $P$ is not $p$--central.

\begin{ex} \label{last group example} Let $P$ be the group of order 64 number \#108 on the Carlson group cohomology website.  From the information there one learns that $Out(P)$ has order $3 \cdot 2^8$, $c(P)=\depth_{H^*(P)}H^*(P)=2$, while $\mrk(P)=\rk(P)=3$.  One can also deduce that $\F_2[Out(P)]$ has precisely two irreducibles - the one dimensional trivial module `1' and a two dimensional module `$S$'.  We compute the following table of nonzero dimensions, where $\emptyset$ denotes that the corresponding summand of $\bar R_d(H^*(P))$ is 0.

\bigskip

\begin{center}
\begin{tabular}{|c|c|c|} \hline
 & & \\[-.12in]
$d$ & $\dim e_1\bar R_dH^*(P)$  & $\dim e_S\bar R_dH^*(P)$ \\   \hline
0 &3 & $\emptyset$ \\   \hline
1 &2 & 3 \\   \hline
2 &2 & 3 \\   \hline
3 &3 & 2 \\   \hline
4 &2 & 2 \\   \hline
5 &2 & 2 \\   \hline
6 &2 & 2 \\   \hline
7 &2 & $\emptyset$ \\   \hline
\end{tabular}
\end{center}

\bigskip

Note that this example shows that, when $d>0$, $e_S\bar R_dH^*(P) \neq 0$ does not imply that $\dim e_S\bar R_dH^*(P) \geq \mrk(P)$, in contrast to \corref{martino priddy cor}.  Indeed, for this group, $\dim \bar R_dH^*(P) = 2 < 3 = \mrk(P)$ for $d=4,5,6,7$.

\end{ex}

\subsection{Organization of the paper} The rest of the paper is organized as follows.  In \secref{previous results section}, we quickly survey previous results related to Conjecture A.  The nilpotent filtration is reviewed in \secref{U section}, which has various relevant general results about $\U$ and $\K$ concerning both dimension and depth. In \secref{HH(P) section}, we give a first proof of \thmref{R_0 theorem} using the functor category description of $\U/Nil$ introduced in \cite{hls1}.  In \secref{Out(P) section}, we then collect various results about $p^{\prime}$--automorphisms of $p$--groups, and prove \thmref{Out(P,V) thm} and \propref{reduced prop}.  In this section we also discuss the family of 2-central groups given in \exref{2C examples}.  In \secref{proof of the main theorem}, we use a convenient formula for $\bar R_dH^*(P)$ from \cite{k4} to prove \thmref{R_d theorem}. A short proof of \propref{pC prop} also appears there.  That Conjecture C implies Conjecture A is shown in the short \secref{conjecture section}, and in \secref{Conj C examples section} various 2--central groups are shown to satisfy Conjecture C, including those listed in \thmref{2 central conj thm}.  Finally, in \secref{examples section}, we will discuss  \exref{last group example} in detail. \\

The author wishes to thank the Cambridge University Pure Mathematics Department for its hospitality during a visit during which a good part of this research was done.  Information compiled by Ryan Higginbottom  has been very useful, and was used by him in \cite{higginbottom} to verify Conjecture A for all groups of order dividing 64 except for group 64\#108 .  \exref{2C examples} arose from a conversation with David Green.

\section{Previous results} \label{previous results section}

In the 1984 paper \cite{ds}, T. Diethelm and U. Stammbach gave a group theoretic proof that, if $P$ is a finite $p$--group, every irreducible $\F_p[Out(P)]$--module occurs as a composition factor in $H^*(P)$.  Independently and concurrently, it was noted explicitly in \cite{hk}, and implicitly in \cite{nishida}, that this same result was a consequence of the proof of the Segal Conjecture in stable homotopy theory.  From both \cite{ds} and \cite{hk}, one can conclude that every irreducible occurs an infinite number of times, so that $\dim e_SH^*(P) \geq 1$ for all irreducibles $S$.

The 1999 paper of P. Symonds \cite{symonds} gives a second group theoretic proof that every irreducible $\F_p[Out(P)]$--module occurs in $H^*(P)$, and inspection of his proof also yields Martino and Priddy's  lower bound $\dim e_SH^*(P) \geq c(P)$.\footnote{\cite{symonds}  was written with knowledge of \cite{hk}, but apparently not of \cite{ds}.}

A proof of Conjecture A for an elementary abelian $p$--group $V$ serves as a starting point for work on the general question, and goes as follows.

Let $S^*(V) \subset H^{2*}(V)$ denote the symmetric algebra generated by $\beta(H^1(V)) \subset H^2(V)$, where $\beta$ is the Bockstein.  Thus $S^*(V)$ is a polynomial algebra on $r= \rk(V)$ generators in degree 2.  It is a classic result of Dickson \cite{dickson, crabb} that the invariant ring $S^*(V)^{GL(V)}$ is again polynomial on homogeneous generators $c_V(1), \dots , c_V(r)$, where $c_V(i)$ has degree $2(p^r-p^{r-i-1})$.

For any idempotent $e \in \F_p[GL(V)]$, $eH^*(V)$ is an $S^*(V)^{GL(V)}$--module via left multiplication.
Now one observes that $S^*(V)$ (and thus $H^*(V)$) is a free $S^*(V)^{GL(V)}$--module, and that every finite $\F_p[GL(V)]$--module occurs as a submodule of $S^*(V)$ \cite[p.45]{alperin}.  These facts imply that, for all irreducible $\F_p[GL(V)]$--modules $S$, $e_SH^*(V)$ is a nonzero finitely generated free $S^*(V)^{GL(V)}$--module, and thus $\dim e_SH^*(V) = \rk(V)$.

Using \cite{hk}, this result immediately extends to all abelian $p$--groups.

In \cite{mp}, Martino and Priddy use the analogous action of $H^*(P)^{Out(P)}$ on $eH^*(P)$ to show that, for all nonzero idempotents $e$, $\dim eH^*(P) \geq \depth_{H^*(P)}H^*(P)$.  Their proof is critically dependent on an unpublished paper that was later withdrawn.  However, as we note in the next section, their argument can be salvaged by using the later work of D. Bourguiba and S. Zarati \cite{bour zarati} utilizing deep properties of $\U$ and $\K$.

\begin{rem} \label{strong conj remark} Readers of \cite{mp} will know that much of Martino and Priddy's paper concerns a refined version of Conjecture A.  The double Burnside ring $A(P,P)$ acts on $H^*(P)$, and one can conjecture that $\dim eH^*(P) = \rk(P)$ for all idempotents $e \in A(P,P)\otimes \F_p$ that project to a nonzero element under the retraction of algebras $A(P,P)\otimes \F_p \ra \F_p[Out(P)]$.  The paper \cite{hk} notes that one can deduce that $eH^*(P)\neq 0$, and thus  that $\dim eH^*(P) \geq 1$, from the Segal conjecture; there is currently no `group theoretic' proof of this fact.  Partly for this reason, we have focused on the $Out(P)$ conjecture, though some of our general theory evidently applies to the $A(P,P)$ version.\footnote{Our results do imply the refined conjecture for $p$--central groups at odd primes.  See \remref{no transfers}.}
\end{rem}

\section{Krull dimension and the nilpotent filtration of $\U$}  \label{U section}

\subsection{The nilclosure functor}

Let $\mathcal Nil_1 \subset \U$ be the localizing subcategory generated by suspensions of unstable $\A$--modules, i.e. $\mathcal Nil_1$ is the smallest full subcategory containing all suspensions of unstable modules that is closed under extensions and filtered colimits.

An unstable module $M$ is called {\em nilreduced} (or just {\em reduced}) if it contains no such nilpotent submodules, or, equivalently, if $\Hom_{\U}(N,M) = 0$ for all $N \in \mathcal Nil_1$.  $M$ is called {\em nilclosed} if also $\Ext^1_{\U}(N,M) = 0$ for all $N \in \mathcal Nil_1$.

Let $L_0: \U \ra \U$ be localization away from $\mathcal Nil_1$.  Thus $L_0M$ is nilclosed, and there is a natural transformation $M \ra L_0M$ with both kernel and cokernel in $\mathcal Nil_1$.

Recall that $T_V: \U \ra \U$ is defined to be left adjoint to $H^*(V) \otimes \text{\underline{\hspace{.1in} }}$.  The various marvelous properties of $T_V$ are reflected in similar properties of $L_0$.

\begin{prop} \label{L0 properties} The functor $L_0: \U \ra \U$ satisfies the following properties. \\

\noindent{\bf (a)} There are natural isomorphisms $ L_0(M \otimes N) \simeq L_0M \otimes L_0N$. \\

\noindent{\bf (b)} There are natural isomorphisms $T_V L_0M \simeq L_0T_V M$. \\

\noindent{\bf (c)} If $K \in \K$, then $L_0K \in \K$, and $K \ra L_0K$ is a map of unstable algebras.
\end{prop}

One approach to properties (a) and (b) is to use \propref{F prop}.  See \cite[I.4.2]{hls2} and \cite{broto zarati 1} for more detail about property (c).

Given a Noetherian unstable algebra $K \in \K$, we recall the category $K_{f.g.}-\U$ as studied in \cite[I.4]{hls2}.  The objects are finitely generated $K$--modules $M$ whose $K$--module structure map $K \otimes M \ra M$ is in $\U$, and morphisms are $K$--module maps in $\U$.

\begin{prop} \label{L0 finiteness} Let $K \in \K$ be Noetherian, and $M \in K_{f.g.}-\U$.  Then $L_0K \in K_{f.g.}-\U$, and thus is Noetherian, and $L_0M \in L_0K_{f.g.}-\U$.
\end{prop}

See \cite[I(4.10)]{hls2}.

\subsection{The nilpotent filtration}

For $d\geq 0$, let $\mathcal Nil_d \subset \U$ be the localizing subcategory generated by $d$--fold suspensions of unstable $\A$--modules.  An unstable module $M$ admits a natural filtration
$$ \dots \subseteq nil_2M \subseteq nil_1M \subseteq nil_0 M= M,$$
where $nil_dM$ is the largest submodule in $\mathcal Nil_d$.

As observed in \cite[Prop.2.2]{k2}, $ nil_dM/nil_{d+1}M = \Sigma^d R_dM$,
where $R_dM$ is a reduced unstable module. (See also \cite[Lemma 6.1.4]{s2}.)

\begin{prop} \label{Rd properties} The functors $R_d: \U \ra \U$ satisfy the following properties. \\

\noindent{\bf (a)} There are a natural isomorphisms $ R_*(M \otimes N) \simeq R_*M \otimes R_*N$ of graded objects in $\U$.  \\

\noindent{\bf (b)} There are natural isomorphism $T_V R_dM \simeq R_dT_V M$.  \\

\noindent{\bf (c)} Let $K \in \K$ be Noetherian, and $M \in K_{f.g.}-\U$.  Then $R_0K$ is also a Noetherian unstable algebra, and $R_dM \in R_0K_{f.g.}-\U$, for all $d$.
\end{prop}

For the first two properties, see \cite[\S 3]{k2}, and the last follows easily from the first.

Now let $\bar R_dM$ denote the nilclosure of $R_d(M)$.  Thus $R_dM \subseteq L_0R_dM = \bar R_dM$.  We also note that $L_0M = \bar R_0M$.

The previous propositions combine to prove the following.

\begin{prop} \label{bar Rd properties} The functors $\bar R_d: \U \ra \U$ satisfy the following properties. \\

\noindent{\bf (a)} There are natural isomorphisms $\bar R_*(M \otimes N) \simeq \bar R_*M \otimes \bar R_*N$ of graded objects in $\U$. \\

\noindent{\bf (b)} There are natural isomorphisms $T_V \bar R_dM \simeq \bar R_dT_V M$.  \\

\noindent{\bf (c)} Let $K \in \K$ be Noetherian, and $M \in K_{f.g.}-\U$.  Then $\bar R_0K$ is also a Noetherian unstable algebra, and $\bar R_dM \in \bar R_0K_{f.g.}-\U$, for all $d$.
\end{prop}

Henn \cite{henn} proved the following important finiteness result.

\begin{prop} \label{finite prop 2} Let $K \in \K$ be Noetherian, and $M \in K_{f.g.}-\U$. Then the nilpotent filtration of $M$ has finite length.  Equivalently, $\bar R_dM = 0$ for $d >> 0$.
\end{prop}

\subsection{Invariant rings as Noetherian unstable algebras}  \label{finite subsection}

% The next result is essentially due to Lannes.

\begin{prop} \label{finite prop 1} Let $K \in \K$ be Noetherian.  Given any subgroup $G < Aut_{\K}(K)$, the invariant ring $K^G$ is again a Noetherian unstable algebra, and $K$ is a finitely generated $K^G$--module.
\end{prop}

\begin{proof}  Let $D(d) \in \K$ denote the $d^{th}$ Dickson algebra: in the notation of \secref{previous results section}, $D(d)$ is the polynomial algebra $S^*((\F_p)^d)^{GL_d(\F_p)}$.  Then let $D(d,j) \in \K$ denote the subalgebra consisting of $p^j$ powers of the elements of $D(d)$.  In \cite[Thm.A.1]{bour zarati}, Lannes shows that, if $d = \dim K$, then there is an embedding of unstable algebras $D(d,j) \ra K$ that is unique in the sense that any two such embeddings will agree after restriction to $D(d,k)$ with $k$ large enough.  As $G$ is necessarily finite, it follows that for large enough $j$, there is an embedding $D(d,j) \ra K^G$ such that $K$ is a finitely generated $D(d,j)$--module.  The proposition follows.
\end{proof}

Our various propositions apply to the case when $K = H^*(P)$.

\begin{cor} Let $P$ be a finite $p$--group. \\

\noindent{\bf (a)} \ Both $H^*(P)^{Out(P)}$ and  $\HH^*(P)^{Out(P)}$ are Noetherian. \\

\noindent{\bf (b)} \  $e_SH^*(P) \in H^*(P)^{Out(P)}_{f.g.}-\U$ and $e_S\bar R_dH^*(P) \in \HH^*(P)^{Out(P)}_{f.g.}-\U$, for all irreducible $\F_p[Out(P)]$--modules $S$ and all $d$.
\end{cor}

\subsection{How to deal with odd primes} \label{odd prime subsection}  When $p$ is odd, $K \in \K$ is not necessarily commutative unless $K$ is concentrated in even degrees.  To use standard definitions and results from the commutative algebra literature, and for other technical reasons, it is useful to have a systematic way of ridding ourselves of this problem.

A standard thing to do, done many times before, and going back at least to \cite{lannes zarati}, goes as follows.  Let $\U^{\prime}$ denote the full subcategory of $\U$ consisting of modules concentrated in even degrees.  Given $M \in \U$, we let $M^{\prime} \in \U^{\prime}$ denote the image of $M$ under the right adjoint of the inclusion $\U^{\prime} \subset \U$: in more down-to-earth terms, $M^{\prime}$ is the largest submodule of $M$ contained in even degrees.  It is easy to see that if $K$ is an unstable algebra, then $K^{\prime}$ is a subalgebra.  As an example, $H^*(V)^{\prime} = S^*(V)$.

Statements about $K_{f.g.}-\U$ become statements about $K^{\prime}_{f.g.}-\U$ using the next result.

\begin{prop} If $K \in \K$ is Noetherian, so is $K^{\prime}$, and $K$ is a finitely generated $K^{\prime}$--module.
\end{prop}

This is \cite[Lemma 5.2]{broto zarati 1}.

\subsection{Krull dimension}

We will generally work with the standard definition of Krull dimension.  Given a commutative graded Noetherian $\F_p$--algebra $K$, $\dim K = d$ if
$$\wp_0 \subset \dots \subset \wp_d$$ is a chain of prime ideals in $K$ of maximal length.  If $M$ is a finitely generated $K$--module, $\dim M$ is then defined to be $\dim (K/Ann(M))$, where $Ann(M)$ is the annihilator ideal of $M$.

Given a Noetherian $K \in \K$ and $M \in K_{f.g.}-\U$, the proposition of the last subsection shows that $M \in K^{\prime}_{f.g.}-\U$, and thus $\dim M$ will be a well defined finite natural number.

In this situation, the Poincar\'e series of $M$, $ \sum_{i=0}^\infty (\dim_{\F_p}M^i) t^i$, will be a rational function, and $\dim M$ equals the order of the pole at $t=1$.  There is a third way of calculating $\dim K$: it is the number $d$ such that there exist algebraically independent elements $k_1, \dots ,k_d \in K$ with $K$ finitely generated over $\F_p[k_1,\dots,k_d]$.  See \cite[\S 10.2]{carlson et al} for a nice discussion of these and related facts.

\propref{Rd prop} of the introduction is a special case of the following.

\begin{prop} \label{general Rd prop} Let $K \in \K$ be Noetherian.  Given $M \in K_{f.g.}-\U$,
$$\dim M = \max_{d} \{ \dim \bar R_dM \}.$$
\end{prop}

\begin{proof}  As the modules $\Sigma^dR_dM$ are the composition factors associated to a finite filtration of $M$, standard properties of Krull dimension \cite[(12.D)]{matsumura} imply that $\displaystyle \dim M = \max_{d} \{ \dim R_dM \}$.  Recalling that $\bar R_d = L_0R_d$, the next proposition finishes the proof.
\end{proof}

\begin{prop} \label{dim L prop} Given $K$ and $M$ as above, if $M$ is reduced then $$\dim M = \dim L_0M.$$
\end{prop}

This follows from the next two lemmas.

\begin{lem} $Ann(M)$ is a sub $\A$--module of $K$.
\end{lem}

\begin{proof}  We are claiming that, if $kx=0$ for all $x \in M$, then, for all $a\in \A$ and $x \in M$, we have $(ak)x=0$.  Fixing $a\in \A$, assume by induction that $(a^{\prime}k)y=0$ for all $y \in M$ and $a^{\prime} \in \A$ with $|a^{\prime}|<|a|$.  With $\Delta a = \sum a^{\prime} \otimes a^{\prime \prime}$, we then have
\begin{equation*}
\begin{split}
0 & = a(kx) \\
  & = (ak)x \ + \sum_{|a^{\prime}|<|a|} (a^{\prime}k)(a^{\prime \prime}x) \\
  & = (ak)x,
\end{split}
\end{equation*}
where the last equality uses the inductive hypothesis.
\end{proof}

\begin{lem} With $M$ reduced, $Ann(M) = Ann(L_0M)$.
\end{lem}

\begin{proof}
\begin{equation*}
\begin{split}
Ann (M) & = \{k\in K \ | \ \A k \otimes M \ra M \text{ is 0} \} \\
  & = \{k\in K \ | \ \A k \otimes M \ra L_0M \text{ is 0} \} \\
  & = \{k\in K \ | \ \A k \otimes L_0M \ra L_0M \text{ is 0} \} \\
  & = Ann (L_0M).
\end{split}
\end{equation*}
Here the first and last equalities are consequences of the last lemma.  The second is immediate since $M$ is reduced, i.e. $M \ra L_0M$ is monic. Finally, the third equality follows from the universal property of nilclosed modules, as $\A k \otimes M \ra \A k \otimes L_0M$ is a $\mathcal Nil_1$--isomorphism.\end{proof}

\subsection{Depth} \label{depth subsection} If $M$ is a $K$--module, and $I \subset K$ is an ideal, the depth of $M$ with respect to $I$ is defined to be the maximal length $l$ of an $M$--regular sequence in $I$: $r_1,\dots, r_l \in I$ such that for each $i$ between 1 and $l$, $r_i$ is not a zero divisor on $M/(r_1,\dots,r_{i-1})M$.

If $K$ is a Noetherian unstable algebra, and $M \in K_{f.g.}-\U$, we let $\depth_KM$ denote the depth of $M$ with respect to the ideal of positive degree elements in $K$.

It is standard that $\depth_KM \leq \dim N$, where $N\subset M$ is any nonzero submodule.  This is usually stated in the following equivalent formulation \cite[Thm.29]{matsumura}: $\depth_KM \leq \dim K/\wp$ where $\wp$ is any associated prime ideal, i.e. a prime ideal arising as the annihilator of an element of $M$. The associated primes include all the minimal primes in the support of $M$ \cite[Thm.9]{matsumura}.

Quillen \cite[Prop.11.2]{quillen part II} shows that the minimal primes correspond to the maximal elementary abelian subgroups, and so $\depth_{H^*(P)}H^*(P) \leq \mrk(P)$, as asserted in the introduction.

The work of Bourguiba, Lannes, and Zarati \cite{bour zarati} shows

\begin{prop}  Let $K \in \K$ be Noetherian, $G < Aut_{\K}(K)$, and $M \in K_{f.g.}-\U$.  Then $\depth_{K^G}M = \depth_KM$.
\end{prop}

The point here is that the main theorem of \cite{bour zarati} says that $\depth_KM$ can be calculated by using a very specific regular sequence of `generalized Dickson invariants', and, as discussed in the proof of \propref{finite prop 1}, Lannes' theorem in the appendix of \cite{bour zarati} says that these elements will be in $K^G$.

\begin{cor}  Let $K \in \K$ be Noetherian, $G < Aut_{\K}(K)$, $M \in K_{f.g.}-\U$, and $N$ be a nonzero direct summand of $M$, viewed as a $K^G$--module.  Then $\dim N \geq \depth_KM$.
\end{cor}
\begin{proof}  $\dim N \geq \depth_{K^G}N \geq \depth_{K^G}M = \depth_KM$.
\end{proof}

Specializing further, we have

\begin{cor}  Let $P$ be a finite $p$--group, and $S$ an irreducible $\F_p[Out(P)]$--module.  If $e_SH^*(P) \neq 0$, then $\dim e_SH^*(P) \geq \depth_{H^*(P)}H^*(P)$.  If $e_S\bar R_dH^*(P) \neq 0$, then   $\dim e_S\bar R_dH^*(P) \geq \depth_{\HH^*(P)}\bar R_d(P)$.
\end{cor}

Duflot showed \cite{duflot} that $\depth_{H^*(P)}H^*(P) \geq \rk(C)=c(P)$.  Subsequent proofs, beginning with \cite{broto henn}, have emphasized that this theorem is a reflection of the $H^*(C)$--comodule structure on $H^*(P)$ induced by the multiplication homomorphism $C \times P \ra P$.  This same homomorphism also induces an $H^*(C)$--comodule structure on $\bar R_*H^*(P)$, and the proof of Duflot's theorem given in \cite[Thm.12.3.3]{carlson et al} goes through without change to prove

\begin{prop} $\depth_{\HH^*(P)}\bar R_*(P) \geq c(P)$.
\end{prop}

This proposition and the preceeding corollary imply \propref{depth prop}.

\section{Proof of \thmref{R_0 theorem}} \label{HH(P) section}

In this section, we give a proof of \thmref{R_0 theorem} based on the identification of $\U/\mathcal Nil$ as a certain functor category.

Following \cite{hls1}, let $\mathcal F$ be the category of functors from finite dimensional $\F_p$--vector spaces to $\F_p$--vector spaces.  This is an abelian category in the obvious way:  $F \ra G \ra H$ is exact if $F(W) \ra G(W) \ra H(W)$ is exact for all $W$.

Let $l: \U \ra \mathcal F$ be defined by $l(M)(W) = (T_WM)^0 = \Hom_{\U}(M,H^*(W))^{\prime}$. (Here $M^{\prime}$ denotes the continuous dual of a profinite vector space $M$.)  This has right adjoint $r: \mathcal F \ra \U$ given by $r(F)^d = \Hom_{\mathcal F}(H_d,F)$, where $H_*(W)$ is the mod $p$ homology of the group $W$.

\begin{prop} \label{F prop} The functors $l$ and $r$ satisfy the following properties. \\

\noindent{\bf (a)} $l$ is exact. \\

\noindent{\bf (b)} The natural transformation $M \ra r(l(M))$ identifies with $M \ra L_0M$. In particular $l(M) = 0$ if and only if $M \in \mathcal Nil_1$.  \\

\noindent{\bf (c)} Both $l$ and $r$ commute with tensor products.

\end{prop}

For proofs of these properties see \cite{hls1, genrep1, s2}.

The proposition implies that a nilclosed module, like $\HH^*(P)$, is completely determined as an object in $\U$ by its associated functor.  It is well known that Quillen's work allows for the identification of $l(H^*(P))= l(\HH^*(P))$.  Let $\Rep(W,P) = \Hom(W,P)/Inn(P)$.

\begin{prop}  $l(\HH^*(P))(W) = \F_p^{\Rep(W,P)}$.
\end{prop}

For a simple proof, see \cite{k1}.

\begin{proof}[Proof of the first statement of \thmref{R_0 theorem}]
By the discussion above, an irreducible $\F_p[Out(P)]$--module $S$ occurs in $\HH^*(P)$ if and only if it appears in the permutation module $\F_p^{\Rep(W,P)}$ for some $W$.

Let $\alpha: W \ra P$ represent an element in $\Rep(W,P)$.  The stablizer of the $Out(P)$--orbit of this element will be $Out_0(P,V)$ where $V = \alpha(W)$.  Thus, as an $Out(P)$--module, $\F_p^{\Rep(W,P)}$ is a direct sum of modules of the form $\F_p^{Out(P)/Out_0(P,V)}$, and every such module appears as a summand of $\F_p^{\Rep(W,P)}$ if $\rk W \geq \rk V$.

Finally, we note that if $V_1 < V_2 < P$, then $Out_0(P,V_2) \leq Out_0(P,V_1)$, and thus $\F_p^{Out(P)/Out_0(P,V_1)}$ is a quotient of $\F_p^{Out(P)/Out_0(P,V_2)}$.  It follows that if $S$ occurs in $\F_p^{Out(P)/Out_0(P,V)}$ for some $V < P$, then $S$ occurs in $\F_p^{Out(P)/Out_0(P,V)}$ for some maximal $V < P$.
\end{proof}

We now finish the proof of \thmref{R_0 theorem} by showing that if $S$ occurs in the permutation module $\F_p^{Out(P)/Out_0(P,V)}$ then $\dim e_S\HH^*(P) \geq \rk V$.  We do this by looking more carefully at $\F_p^{\Rep(W,P)}$ as an object in $\mathcal F$ with $Out(P)$ action.

We begin by giving an increasing filtration of the set valued bifunctor $\Rep(W,P)$.  Let $\Rep_k(W,P) \subseteq \Rep(W,P)$ be the set of elements represented by homomorphisms with image of rank at most $k$.  The inclusions
$$ \{0\} = \Rep_0(W,P) \subseteq \Rep_1(W,P) \subseteq \dots \subseteq \Rep_{\rk(P)}(W,P)= \Rep(W,P)$$
induces epimorphisms of objects in $\mathcal F$ with $Out(P)$ action
$$ \F_p = \F_p^{\Rep_0(W,P)} \la  \F_p^{\Rep_1(W,P)} \la \cdots \la  \F_p^{\Rep_{\rk(P)}(W,P)}=  \F_p^{\Rep(W,P)}.$$

An appeal to \propref{dim L prop} thus shows

\begin{lem} $\displaystyle \dim e_S\HH^*(P) = \max_k \{ \dim e_S r (\ker \{\F_p^{\Rep_k(W,P)} \ra \F_p^{\Rep_{k-1}(W,P)} \}) \}$.
\end{lem}

To more usefully describe $\ker \{\F_p^{\Rep_k(W,P)} \ra \F_p^{\Rep_{k-1}(W,P)} \}$, we need to introduce some notation.  Given $V < P$, let $GL_P(V) = Out(P,V)/Out_0(P,V)$.  Note that $GL_P(V)$ is naturally a subgroup of $GL(V)$ and also acts on the right of the $Out(P)$--set $Out(P)/Out_0(P,V)$.

\begin{lem}  There is an isomorphism of objects in $\mathcal F$ with $Out(P)$ action,
$$ \ker \{\F_p^{\Rep_k(W,P)} \ra \F_p^{\Rep_{k-1}(W,P)} \} \simeq \prod_{[V]} \F_p^{\Epi(W,V) \times_{GL_P(V)} Out(P)/Out_0(P,V)},$$
where the product is over $Aut(P)$--orbits of elementary abelian subgroups of rank $k$.
\end{lem}

Now we need to identify the nilclosed unstable module with $Out(P)$ action associated to $\F_p^{\Epi(W,V) \times_{GL_P(V)} Out(P)/Out_0(P,V)}$.

Let $c_V \in S^*(V)^{GL(V)}$ be the `top' Dickson invariant: $c_V$ is the product of all the elements $\beta(x)\in H^2(V)$, with $0 \neq x \in H^1(V)$.  The key property we need is that $\eta^*(c_V) = 0$ for all proper inclusions $\eta: U < V$.

\begin{lem}[Compare with {\cite[proof of Thm.II.6.4]{hls1}}] $l(c_VS^*(V))(W) = \F_p^{\Epi(W,V)}$.
\end{lem}

\begin{cor} The functor $l$ assigns to the reduced module
$$[c_VS^*(V) \otimes \F_p^{Out(P)/Out_0(P,V)}]^{GL_P(V)},$$
the object in $\mathcal F$ which sends $W$ to $\F_p^{\Epi(W,V) \times_{GL_P(V)} Out(P)/Out_0(P,V)}.$
\end{cor}

To finish the proof of \thmref{R_0 theorem}, we are left needing to prove

\begin{prop} If an irreducible $\F_p[Out(P)]$--module $S$ occurs in the permutation module $\F_p[Out(P)/Out_0(P,V)]$, then
$$\dim e_S [c_VS^*(V) \otimes \F_p^{Out(P)/Out_0(P,V)}]^{GL_P(V)} = \rk V.$$
\end{prop}

\begin{proof}  The proof is similar to the proof of Conjecture A when $P=V$ given in \secref{previous results section}.  Our first observation is that $e_S [c_VS^*(V) \otimes \F_p^{Out(P)/Out_0(P,V)}]^{GL_P(V)}$ is a $S^*(V)^{GL(V)}$--submodule of the free $S^*(V)^{GL(V)}$--module $S^*(V) \otimes \F_p^{Out(P)/Out_0(P,V)}$.  As $\dim S^*(V)^{GL(V)} = \rk V$, it suffices to show that if $S$ occurs in the permutation module $\F_p[Out(P)/Out_0(P,V)]$, then $S$ occurs in
$$[c_VS^*(V) \otimes \F_p^{Out(P)/Out_0(P,V)}]^{GL_P(V)}.$$
The group ring $\F_p[GL_P(V)]$ occurs as a submodule of $S^*(V)$ \cite[p.45]{alperin}, and thus as a submodule of $c_VS^*(V)$. Thus $S$ will occur in $[c_VS^*(V) \otimes \F_p^{Out(P)/Out_0(P,V)}]^{GL_P(V)}$ if it occurs in
$[\F_p[GL_P(V)] \otimes \F_p^{Out(P)/Out_0(P,V)}]^{GL_P(V)}$.  But this last $\F_p[Out(P)]$--module rewrites as $\F_p^{Out(P)/Out_0(P,V)}$.
\end{proof}

\section{Proof of \thmref{Out(P,V) thm} and related results}     \label{Out(P) section}

\thmref{R_0 theorem} tells us that calculating $\dim e_S \HH^*(P)$ amounts to understanding the $\F_p[Out(P)]$--module composition factors of the permutation modules
$$\F_p[Out(P)/Out_0(P,V)],$$ when $V < P$ is a maximal elementary abelian subgroup.

In this section, we use group theory to find various conditions ensuring that all irreducibles occur in this way.

\subsection{Composition factors of permutation modules}
We begin with an elementary, but useful, lemma about permutation modules.

\begin{lem}  \label{perm lemma} Let $H$ be a subgroup of a finite group $G$. Then {\bf (a) $\Rightarrow$ (b) $\Rightarrow$ (c)}, and, if $H$ is also normal, {\bf (c) $\Rightarrow$ (a)}. \\

\noindent{\bf (a)} $H$ is a $p$--group. \\

\noindent{\bf (b)} Every irreducible $\F_p[G]$--module occurs as a submodule of $\F_p[G/H]$. \\

\noindent{\bf (c)} Every irreducible $\F_p[G]$--module occurs as a composition factor of $\F_p[G/H]$.
\end{lem}
\begin{proof} We begin by reminding the reader that permutation modules are self dual, so an irreducible $S$ occurs as a submodule of $\F_p[G/H]$ if and only if $S$ occurs as a quotient module.  This happens if and only if there exists a nonzero homomorphism $\F_p[G/H] \ra S$, which is equivalent to $S^H \neq 0$.

If $H$ is a $p$--group, then $M^H \neq 0$ for all nonzero $\F_p[H]$--modules $M$.  Thus (a) implies (b), and (b) implies (c) is obvious.

Now suppose that $H$ is normal.  Then the composition factors of $\F_p[G/H]$ are just the irreducible $\F_p[G/H]$--modules pulled back to $G$.  But if $H$ is not a $p$--group, then $G/H$ has fewer $p^{\prime}$ conjugacy classes than $G$, and thus fewer irreducible modules.  Thus (c) implies (a) under the normality assumption.
\end{proof}

\begin{rem}  It is not hard to see that, if $H$ is not normal, then neither implication (b) $\Rightarrow$ (a) nor (c) $\Rightarrow$ (b) need hold.  For the former, let $p=3$, $G = SL_2(\F_3)$, and $H = \Z/2$ permuting a basis for $(\F_3)^2$.  For the latter, let $p=3$, $G=\Sigma_3$, and $H = \Z/2$.
\end{rem}

\begin{cor} \label{little cor} If a $p$--group $P$ has a unique maximal elementary abelian subgroup $V$, then every irreducible occurs in $\F_p[Out(P)/Out_0(P,V)]$ if and only if $Out_0(V)$ is a $p$--group.
\end{cor}

\begin{proof} The hypothesis implies that $Out_0(P,V)$ arises as the kernel of an evident homomorphism $Out(P) \ra Aut(V)$, and thus is normal in $Out(P)$. Now the lemma applies.
\end{proof}

\subsection{$p^{\prime}$--automorphisms of $p$--groups}  \label{p' subsection}

We will use various well known results about detecting automorphisms of a $p$--group $P$ of  order prime to $p$.

The most classic is due to Burnside and Hall, and goes as follows.  Let $\Phi(P)$ be the Frattini subgroup of $P$, so that $P/\Phi(P) = H_1(P;\F_p)$.

\begin{prop}[{\cite[Thm.5.1.4]{gor}}] \label{phi(P) prop} The kernel of the homomorphism $Aut(P) \ra Aut(P/\Phi(P))$ is a $p$--group.
\end{prop}

With groups $A$ and $B$ as indicated, the Thompson $A \times B$ lemma (\cite[(24.2)]{asch}, \cite[Thm.5.3.4]{gor}) immediately applies to prove the next proposition.

\begin{prop} \label{CP(B) lemma} Let $B$ be an abelian subgroup of a finite $p$--group $P$, and let $A = \{ \alpha \in Aut(P) \ | \ \alpha(g) = g \text{ for all } g\in C_P(B)\}$.  Then $A$ is a $p$--group.
\end{prop}

\begin{cor}  \label{CP(V) prop} For all elementary abelian $V < P$, the kernel of the homomorphism
$ Aut(P,V) \ra Aut(C_P(V))$ is a $p$--group.  Thus if $Aut(C_P(V))$ is a $p$--group so is $Aut(P,V)$, and if $Aut_0(C_P(V),V)$ is a $p$--group so is $Aut_0(P,V)$.
\end{cor}

The next result requires that  $p$ be odd.  Let $\Omega_1(P)$ be the subgroup of $P$ generated by the elements of order $p$.

\begin{prop}[{\cite[Thm.5.3.10]{gor}}] \label{omega prop}  Let $P$ be a $p$--group, with $p$ odd.  Then the kernel of the homomorphism $ Aut(P) \ra Aut(\Omega_1(P))$ is a $p$--group.
\end{prop}

\begin{cor} \label{pC Aut cor} Let $Q$ is $p$--central, with $p$ odd, and with maximal central elementary abelian subgroup $C$.  Then $Aut_0(Q,C)$ is a $p$--group.
\end{cor}

This follows from the proposition, noting that $C = \Omega_1(Q)$, since $Q$ is $p$--central.

\begin{proof}[Proof of \thmref{Out(P,V) thm}] The first statement follows immediately from \corref{CP(V) prop}. Now we note that, if $V < P$ is a maximal elementary abelian subgroup, then $C_P(V)$ is $p$--central with $V$ as its maximal central elementary abelian subgroup.  Thus \corref{pC Aut cor} applies to prove the second statement.
\end{proof}

\subsection{$2$--central group examples} \label{2C example section}

We elaborate on \exref{2C examples}.

\begin{ex} Let $G_t$ be the $2$--Sylow subgroup of $SU_3(\F_{2^t})$.  This is a $2$--central group of rank $t$, and if $C_t$ is its maximal elementary abelian subgroup, $Out_0(G_t,C_t)$ contains a cyclic subgroup of order $2^t+1$, and thus is not a 2--group.

This will be a consequence of a few facts\footnote{The authors are very grateful to David Green for showing us his unpublished `Appendix B' to \cite{green}, which is our source of information about these groups.} about $G_t$.

Let $q=2^t$. Given $a \in \F_{q^2}$, let $\bar a = a^q$, so that $\F_q$ identifies with the set of $b \in \F_{q^2}$ such that $b+\bar b = 0$.   Then $G_t \subset SU_3(\F_{q}) \subset GL_3(\F_{q^2})$ can be described as
$$ G_t = \left\{
A(a,b) = \left(
\begin{array}{ccc}
1 & a & b \\
0 & 1 & \bar a \\
0 & 0 & 1 \end{array}
\right) \ : \  a,b \in \F_{q^2} \text{ with } b+\bar b = a\bar a \right\}.
$$
The center $C_t$ of $G_t$ is the set of  matrices of the form $A(0,b)$: note that then $b \in \F_{q}$.  All other elements have order 4, and thus  $G_t$ is 2--central sitting in the central extension
$$ C_t \ra G_t \ra G_t/C_t,$$
with $C_t \simeq \F_q$ and $G_t/C_t \simeq \F_{q^2}$.

The normalizer of $G_t$ in $SU_3(\F_{2^t})$ is the semidirect product $G_t >\hspace{-.08in}\lhd T_t$, where $T_t$ is the set of matrices
$$ \left\{
D(c) = \left(
\begin{array}{ccc}
c & 0 & 0 \\
0 & \bar c c^{-1} & 0 \\
0 & 0 & \bar c^{-1} \end{array}
\right) \ : \  c \in \F_{q^2}^{\times} \right\}.
$$
Direct computation shows that
$$D(c)A(a,b)D(c^{-1}) = A({\bar c}^{-1}c^2a,\bar ccb) = A(c^{2-q}a,c^{q+1}b).$$
From this, one deduces that $T_t \ra Aut(G_t/C_t)$, and thus $T_t \ra Aut(G_t)$, is monic if $\gcd(q-2,q^2-1) =1$. This is the case unless $q = 2$, i.e., $t =1$.  Meanwhile, the kernel of the homomorphism $T_t \ra Aut(C_t)$ identifies with the kernel of the multiplicative norm $\F_{q^2}^{\times} \ra \F_q^{\times}$, and thus is cyclic of order $q+1$.  Thus, if $t \geq 2$, $Aut_0(G_t,C_t) = \ker \{ Aut(G_t) \ra Aut(C_t)\}$ contains a cyclic group of order $(q+1)$.

(When $t=1$, $G_t = Q_8$, and so $Aut_0(G_t,C_t)$ also contains a group of order $(q+1)=3$.)
\end{ex}

\begin{rem}  The group $G_2$ has been of interest to those studying group cohomology.  Even though it is $2$--central, a presentation of its cohomology ring is remarkably nasty to write down: see the calculations for group number \# 187 of order 64 in \cite{carlson website, carlson et al}.  It is the smallest group with nontrivial products in its  essential cohomology \cite{green}.  The calculation of the nilpotent length of $p$--central groups in \cite{k3,k4} shows that its nilpotent length is 14, which seems likely to be maximal among all groups of order 64.  Presumably, the groups $G_t$ for larger $t$ are similarly interesting.
\end{rem}

\subsection{$p$--groups with reduced cohomology}

Recall that $P$ has reduced cohomology if and only if $H^*(P)$ is detected by restriction to the various $V <P$, i.e., $H^*(P) \ra \HH^*(P)$ is monic.

We restate \propref{reduced prop}.

\begin{prop}  If $H^*(P)$ is reduced, and $V < P$ is a maximal elementary abelian subgroup, then $C_P(V) = V$.
\end{prop}

Recalling that $C_P(V)$ is $p$--central if $V<P$ is maximal, the proposition is a consequence of the next two lemmas.

\begin{lem} \label{lemma 1} If $H^*(P)$ is reduced, so is $H^*(C_P(V))$ for any elementary abelian $V < P$.
\end{lem}

\begin{lem} \label{lemma 2} If $Q$ is a $p$--central $p$--group with reduced cohomology, then $Q$ is elementary abelian.
\end{lem}

\begin{proof}[Proof of \lemref{lemma 1}]  Properties of the functor $T_V$ imply that if $M \in \U$ is reduced so is $T_VM$: recalling that $T_V$ preserves monomorphisms, property (b) of \propref{L0 properties} impies this. Now one uses that, given $V < P$, $H^*(C_P(V))$ is a (canonical) direct summand in $T_VH^*(P)$. (See \cite[I(5.2)]{hls2} for a calculation of $T_VH^*(P)$.)
\end{proof}

\begin{proof}[Proof of \lemref{lemma 2}]  Let $C$ be the center of $Q$.  The hypothesis is that $H^*(Q) \ra H^*(C)$ is monic, which implies that the inflation map $H^*(Q/C) \ra H^*(Q)$ is zero.  But $H^1(Q/C) \ra H^1(Q)$ is always monic, so we conclude that $H^1(Q/C) = 0$.  But, since $Q/C$ is a $p$--group, this means that $Q/C$ is trivial, i.e., $C = Q$.
\end{proof}

\subsection{A remark about odd prime $p$--central groups} \label{no transfers}

Using a result of Henn and Priddy, \thmref{Out(P,V) thm} has an addendum. For odd prime $p$--central groups $P$, the refinement of Conjecture A described in \remref{strong conj remark} holds: indeed, $\dim e\HH^*(P) = \rk(P)$ for all idempotents $e \in A(P,P)\otimes \F_p$ that project to a nonzero element under the retraction of algebras $A(P,P)\otimes \F_p \ra \F_p[Out(P)]$.  Here dimension is defined via the Poincar\'e series, and we remind the reader that $A(P,P)$ is the double Burnside ring.

The argument goes as follows.  Given an irreducible $\F_p[Out]$--module $S$, let $e_S \in \F_p[Out(P)]$ and $\tilde e_S \in A(P,P)\otimes \F_p$ be the associated primitive idempotents\footnote{Thus $\tilde e_S$ projects to $e_S$, and, in the terminology of \cite{mp}, $e_SBP$ is the super dominant summand associated to the dominant summand $\tilde e_SBP$.}.  Then $e_S\HH^*(P) = \tilde e_S\HH^*(P) \oplus M^*$ where $M^*$ is a finite direct sum of modules of the form $\tilde e_T\HH*(Q)$ for appropriate pairs $(Q,T)$.  By \cite[Prop.1.6.1 and  Lem.2.1]{henn priddy}, if $P$ is an odd prime $p$--central group, the only $Q$'s that can occur here will be proper retracts of $P$.  But such groups will necessarily have strictly smaller ranks.  Thus $\dim M^* < \rk(P)$, and so $\dim \tilde e_S\HH^*(P) = \dim e_S\HH^*(P) = \rk(P)$.

\section{Proof of \thmref{R_d theorem}} \label{proof of the main theorem}

\subsection{A formula for $\bar R_dH^*(P)$}

Recall that $\bar R_0H^*(P) = \lim_{V \in \A^C(P)}H^*(V).$  Using the fact that morphisms in $\A(V)$ factor as inner automorphisms composed with inclusions, this rewrites as follows: there is a natural isomophism
$$ \bar R_0H^*(P)  = \text{Eq} \left\{ \left[\prod_{V} H^*(V)\right]^{Inn(P)}
\begin{array}{c} \mu \\[-.08in] \longrightarrow \\[-.1in] \longrightarrow \\[-.1in] \nu
\end{array}
\prod_{V_1 < V_2} H^*(V_1) \right\},
$$
where the products are over objects and inclusions in $\A^C(P)$, and $\mu$ and $\nu$ are induced by
$ 1: H^*(V_1)\ra H^*(V_1)$ and $ \eta^*: H^*(V_2) \ra H^*(V_1)$
for each inclusion $\eta: V_1 < V_2$ in $\A^C(P)$.

We begin this section with a formula for $\bar R_dH^*(P)$ that generalizes this.

The ingredients of our formula are the following.

First, recall our notation from the introduction: if $W$ is a central elementary abelian $p$--subgroup of $Q$, then $P_WH^*(Q)$ denotes the primitives in the $H^*(W)$--comodule $H^*(Q)$.

An automorphism $\alpha: P \ra P$ induces an isomorphism of unstable algebras $P_{\alpha(V)}H^*(C_P(\alpha(V))) \ra P_VH^*(C_P(V))$.  Thus $P_VH^*(C_P(V))$ is an $Aut(P,V)$--module, and the full automorphism group $Aut(P)$ acts on the product $$\displaystyle \prod_{V \in \A^C(P)} H^*(V) \otimes P_VH^d(C_P(V)).$$

Similarly, an inclusion $\eta: V_1 < V_2$ induces an inclusion $C_P(V_2) < C_P(V_1)$, and thus a map of unstable algebras $\eta_*: P_{V_1}H^*(C_P(V_1)) \ra P_{V_1}H^*(C_P(V_2))$.  We will also denote by $i: P_{V_2}H^*(C_P(V_2)) \ra P_{V_1}H^*(C_P(V_2))$ the evident inclusion of primitives.

Starting from formulae in \cite{hls2}, in \cite{k4}, the author will deduce the following formula for $\bar R_dH^*(P)$.

\begin{prop} \label{good Rd formula prop} $\bar R_dH^*(P)$ is naturally isomorphic to the equalizer
$$ \text{Eq} \left\{ \left[\prod_{V} H^*(V) \otimes P_VH^d(C_P(V))\right]^{Inn(P)}
\begin{array}{c} \mu \\[-.08in] \longrightarrow \\[-.1in] \longrightarrow \\[-.1in] \nu
\end{array}
\prod_{V_1 < V_2} H^*(V_1) \otimes P_{V_1}H^d(C_P(V_2)) \right\},
$$
where the products are over objects and proper inclusions in $\A^C(P)$, and $\mu$ and $\nu$ are induced by
$$ 1 \otimes \eta_*: H^*(V_1) \otimes P_{V_1}H^d(C_P(V_1)) \ra H^*(V_1) \otimes P_{V_1}H^d(C_P(V_2))$$
and
$$ \eta^* \otimes i: H^*(V_2) \otimes P_{V_2}H^d(C_P(V_2)) \ra H^*(V_1) \otimes P_{V_1}H^d(C_P(V_2))$$
for each proper inclusion $\eta: V_1 < V_2$ in $\A^C(P)$.
\end{prop}

Otherwise said, $ x = (x_V) \in [\prod_V H^*(V) \otimes P_VH^d(C_P(V))]^{Inn(P)}$ is in $\bar R_dH^*(P)$ exactly when the components are related by $(1 \otimes \eta_*)(x_{V_1}) = (\eta^* \otimes i)(x_{V_2})$, for each proper inclusion $\eta: V_1 < V_2$ in $\A^C(P)$.

\subsection{ $\bar R_d H^*(Q)$ when $Q$ is $p$--central}

The formula for $\bar R_dH^*(Q)$ simplies dramatically if $Q$ is $p$--central, as the first product in the equalizer formula is indexed by one element, and the second is empty.

\begin{prop}  Suppose $Q$ is $p$--central with maximal central elementary abelian subgroup $C$.  Then
$$ \bar R_dH^*(Q) \simeq H^*(C) \otimes P_CH^d(Q).$$
\end{prop}

As a consequence we can prove \propref{pC prop}, which said that if $P$ is 2--central, then every irreducible $\F_2[Out_0(Q,C)]$--module occurs as a composition factor in $P_CH^*(Q)$.

\begin{proof}[Proof of \propref{pC prop}]  Every irreducible $\F_2[Out(Q)]$--module occurs in $H^*(Q)$, so the same is true for irreducible $\F_2[Out_0(Q,C)]$--modules.  Thus every irreducible $\F_2[Out_0(Q,C)]$--module occurs in $\bar R_*H^*(Q) \simeq H^*(C) \otimes P_CH^*(Q)$.  But, by definition, $H^*(C)$ is a trivial $\F_2[Out_0(Q,C)]$--module, so the irreducibles must occur in the other factor, $P_CH^*(Q)$.
\end{proof}

\begin{rem} In \cite{k4}, we will prove that if $Q$ is $p$--central, the finite unstable algebra $P_CH^*(Q) \subset H^*(Q)$ identifies with the subalgebra of locally finite elements: $x \in H^*(Q)$ such that $\A x \subset H^*(Q)$ is finite.
\end{rem}

\subsection{Proof of \thmref{R_d theorem}}
It will be useful to let $M^*(V)$ denote $P_VH^*(C_P(V))$. Similarly, given $V_1 < V_2$, let $M^*(V_2,V_1) = P_{V_1}H^*(C_P(V_2))$.

Given $V < P$, define $\rho_V$ to be the composite
$$ \bar R_d H^*(P) \subseteq [\prod_{U} H^*(U) \otimes M^d(U)]^{Inn(P)} \lra [H^*(V) \otimes M^d(V)]^{W_P(V)},$$
where the first map is the inclusion of $\F_p[Out(P)]$--modules given in \propref{good Rd formula prop}, and the second is the canonical projection of $\F_p[Out(P,V)]$--modules.  Then let
$$ \tilde \rho_V: \bar R_d H^*(P) \ra \Ind_{Out(P,V)}^{Out(P)}([H^*(V) \otimes M^d(V)]^{W_P(V)})$$
to be the $\F_p[Out(P)]$--module map adjoint to the $\F_p[Out(P,V)]$--module map $\rho_V$.\footnote{Recall that if $H$ is a subgroup of a finite group $G$, induction $\Ind_H^G$ is both left and right adjoint to restriction $\Res_H^G$.}

\begin{defn}  Given an irreducible $\F_p[Out(P)]$--module $S$, and $d \geq 0$, define
$$ {\sup}_d(S) = \{ V \in \text{ob } \A^C(P) \ | \ S \text{ occurs in } \tilde \rho_V(\bar R_d H^*(P)) \}.$$
\end{defn}

It is easy to check that $V \in \sup_d(S)$ if and only if $\rho_{\alpha(V)}(e_S\bar R_d H^*(P)) \neq 0$ for some $\alpha \in Aut(P)$.  Equivalently, $V \in \sup_d(S)$ if and only if there exists $x = (x_U) \in e_S \bar R_d H^*(P)$ with $x_{\alpha(V)} \neq 0$ for some $\alpha \in Aut(P)$.

We will show that $\sup_d(S)$ has the properties described in \thmref{R_d theorem}.

Clearly, given $\alpha \in Aut(P)$, $V \in \sup_d(S)$ implies that $\alpha(V) \in \sup_d(S)$.  Thus $\sup_d(S)$ is a union of $Aut(P)$ orbits in $\ob \A^C(P)$, i.e., property (a) of \thmref{R_d theorem} holds.

Also evident is the first statement of property (c): $\sup_d(S) \neq \emptyset$ if and only if $e_S\bar R_dH^*(P) \neq 0$.

Not much harder is the proof of property (b), which we restate here as a lemma.

\begin{lem}  \label{prop b lemma} Given $\eta: V_1 < V_2$, suppose that
$\eta_*: M^d(V_1) \ra M^d(V_2,V_1)$
is monic.  Then $V_1 \in \sup_d(S)$ implies that $V_2 \in \sup_d(S)$.
\end{lem}

\begin{proof} Suppose $x = (x_U) \in e_S \bar R_d H^*(P)$ has $x_{\alpha(V_1)} \neq 0$, with $\alpha \in Aut(P)$.  Let $\alpha(\eta): \alpha(V_1) < \alpha(V_2)$ be the inclusion.  Then $\eta_*$ monic implies that $\alpha(\eta)_*$ is also monic. Since $(\alpha(\eta)^* \otimes i)(x_{\alpha(V_2)}) = (1 \otimes \alpha(\eta)_*)(x_{\alpha(V_1)}))$, it follows that $x_{\alpha(V_2)} \neq 0$.
\end{proof}

Now suppose $e_S\bar R_dH^*(P) \neq 0$.  As $e_S\bar R_dH^*(P)$ embeds in $\prod_V  \rho_V( e_S\bar R_dH^*(P))$,
$$\dim e_S\bar R_dH^*(P) = \max \{ \dim \rho_V(e_S\bar R_dH^*(P)) \ | \ V \in {\sup}_d(S) \}.$$
Furthermore, $\dim \tilde \rho_V(e_S\bar R_dH^*(P)) = \max \{ \dim \rho_{\alpha(V)}(e_S\bar R_dH^*(P)) \ | \ \alpha \in Aut(P) \},$ and there is an obvious bound: $\dim \rho_V(e_S\bar R_dH^*(P)) \leq \dim H^*(V) = \rk(V)$.
Thus the next two lemmas prove properties (c) and (d) of \thmref{R_d theorem}.

\begin{lem} \label{prop c lemma}  If $V$ is maximal in $\sup_d(S)$, then $\dim \tilde \rho_V(e_S\bar R_dH^*(P)) \geq \rk(V)$.
\end{lem}

\begin{lem} \label{prop d lemma}  If $V$ is maximal in $\ob \A^C(P)$, then $S$ occurs in $\tilde \rho_V(\bar R_dH^*(P))$ if and only if $S$ occurs in $\displaystyle \Ind_{Out(P,V)}^{Out(P)}([\F_p[Aut(P,V)/Aut_0(P,V)] \otimes M^d(V)]^{W_P(V)})$.
\end{lem}

To prove these, we need to introduce yet more notation.

Given $V < P$, let $\Ess^*(V) = \ker \{ M^*(V) \ra \prod_{V<U} M^*(U) \}$, where the product is over all proper inclusions.  Then $\alpha_*(\Ess^*(V)) = \Ess^*(\alpha(V))$ for all $\alpha \in Aut(P)$.  In particular, $\Ess^*(V)$ is an $\F_p[Aut(P,V)]$--submodule of $M^*(V)$.

\begin{lem} \label{factor lemma} There is a map of $Out(P)$--modules
$$\tilde f_V: \Ind_{Out(P,V)}^{Out(P)}([c_VH^*(V) \otimes \Ess^d(V)]^{W_P(V)}) \ra \bar R_dH^*(P)$$
such that $\tilde \rho_V \circ \tilde f_V$ is the canonical inclusion induced by
$$c_VH^*(V) \otimes \Ess^d(V) \subseteq H^*(V) \otimes M^*(V).$$
\end{lem}

\begin{proof}  It is equivalent to define an $Out(P,V)$--module map
$$f_V: [c_VH^*(V) \otimes \Ess^d(V)]^{W_P(V)} \ra \bar R_dH^*(P)$$
such that $\tilde \rho_V \circ f_V$ is an appropriate inclusion.  Given $x \in [c_VH^*(V) \otimes \Ess^d(V)]^{W_P(V)}$ define $f_V(x) \in \left[\prod_{U} H^*(U) \otimes M^d(U) \right]^{Inn(P)}$
by letting its $U^{th}$ component, $f_V(x)_U$, be $({\alpha^{-1*}} \otimes \alpha_*)(x)$ if $U = \alpha(V)$ with $\alpha \in Inn(P)$, and 0 otherwise.

We claim that $f_V(x)$ is, in fact, in $\bar R_dH^*(P)$.  For given $\eta: V_1 < V_2$ with $V_1$ conjugate to $V$, $(1 \otimes \eta_*)(f_V(x)_{V_1}) = 0$ because $\eta_*(Ess^d(V)) = 0$.  Similarly, given $\eta: V_1 < V_2$ with $V_2$ conjugate to $V$, $(\eta^* \otimes i)(f_V(x)_{V_2}) = 0$ because $\eta^*(c_V) = 0$.

Finally, checking that $\tilde \rho_V \circ f_V$ is what it should be is straightforward.
\end{proof}

Next we observe that $\Ind_{Out(P,V)}^{Out(P)}([H^*(V) \otimes M^d(V)]^{W_P(V)})$ is a $H^*(V)^{GL(V)}$--module with the property that each nonzero cyclic $S^*(V)^{GL(V)}$--submodule is free.

\begin{proof}[Proof of \lemref{prop c lemma}]  If $V$ is maximal in $\sup_d(S)$, then, after possibly replacing $V$ by $\alpha(V)$ for some $\alpha \in Aut(P)$, there exists $x \in e_S\bar R_sH^*(V)$ with $x_V \neq 0$ in $[H^*(V) \otimes M^d(V)]^{W_P(V)}$, and $V$ maximal with this property.  Reasoning as in the proof of \lemref{prop b lemma}, it follows that $x_V \in [H^*(V) \otimes \Ess^d(V)]^{W_P(V)}$.
By the proof of \lemref{factor lemma},
$$  c_VS^*(V)^{GL(V)}\tilde \rho_V(x) \subset
\tilde \rho_V(\bar R_dH^*(P)). $$
By the observation above, we also have
$$ c_VS^*(V)^{GL(V)}\tilde \rho_V(x) \subset e_S\Ind_{Out(P,V)}^{Out(P)}([H^*(V) \otimes M^d(V)]^{W_P(V)}).$$
Thus $ c_VS^*(V)^{GL(V)}\tilde \rho_V(x)$ is contained in the intersection of these, i.e.
$$ c_VS^*(V)^{GL(V)}\tilde \rho_V(x) \subset \tilde \rho_V(e_S\bar R_dH^*(P)).$$
But, as a graded $\F_p$--vector space, $c_VS^*(V)^{GL(V)}\tilde \rho_V(x)$ is just the $(|c_V| + |x|)^{th}$ suspension of $S^*(V)^{GL(V)}$, and so $$\rk(V) = \dim c_VS^*(V)^{GL(V)}\tilde \rho_V(x) \leq  \dim \tilde \rho_V(e_S\bar R_dH^*(P)).$$
\end{proof}

\begin{proof}[Proof of \lemref{prop d lemma}]  If $V$ is maximal in $\A^C(P)$, then the product defining $\Ess^*(V)$ is empty, so $\Ess^*(V) = M^*(V)$.  From \lemref{factor lemma}, we thus deduce that we have inclusions of $\F_p[Out(P)]$--modules
$$ \Ind_{Out(P,V)}^{Out(P)}([c_VH^*(V) \otimes M^d(V)]^{W_P(V)}) \subset \tilde \rho_V(\bar R_dH^*(P))$$
and
$$ \tilde \rho_V(\bar R_dH^*(P)) \subset \Ind_{Out(P,V)}^{Out(P)}([H^*(V) \otimes M^d(V)]^{W_P(V)}).$$
Thus the next lemma will finish the proof.
\end{proof}

\begin{lem} \label{occurs lemma} Let $S$ be an irreducible $\F_p[Out(P)]$--module, and $M$ be an $\F_p[Aut(P,V)]$--module.  Then the following are equivalent. \\

\noindent{\bf (a)} $S$ occurs in $\Ind_{Out(P,V)}^{Out(P)}([c_VH^*(V) \otimes M]^{W_P(V)})$. \\

\noindent{\bf (b)} $S$ occurs in $\Ind_{Out(P,V)}^{Out(P)}([H^*(V) \otimes M]^{W_P(V)})$. \\

\noindent{\bf (c)} $S$ occurs in $\Ind_{Out(P,V)}^{Out(P)}([\F_p[Aut(P,V)/Aut_0(P,V)] \otimes M]^{W_P(V)})$.
\end{lem}

\begin{proof}  Noting once again that multiplication by $c_V$ is a monomorphism on the module in (b), we see that (a) and (b) are equivalent.

To simplify notation, let $G = Aut(P,V)/Aut_0(P,V)$ and $W = W_P(V)$.  Thus $W < G < GL(V)$.  The equivalence of (b) and (c) follows immediately from the following claim: if $T$ is an irreducible $\F_p[G/W]$--module, then $T$ occurs in $(H^*(V) \otimes M)^W$ if and only if $T$ occurs in $(\F_p[G] \otimes M)^W$.  The claim is a consequence of the well known fact that, for any $G < GL(V)$, every finite $\F_p[G]$--module embeds in $S^*(V)$ \cite{alperin}, and thus also in $H^*(V)$.  Thus $T$ occurs in $(H^*(V) \otimes M)^W$ if and only if $T$ occurs in $(N \otimes M)^W$ for some finite $\F_p[G]$--module $N$, and this happens if and only if $T$ occurs in (a direct sum of copies of) the module $(\F_p[G] \otimes M)^W$.
\end{proof}

\section{Conjecture C implies Conjecture A} \label{conjecture section}

If $H$ is a subgroup of $G$, and $M$ is an $\F_p[H]$--module, it is easy to verify the following: if every irreducible $\F_p[H]$--module occurs in $M$ then every irreducible $\F_p[G]$--module occurs in $\Ind_H^GM$.
In light of \corref{conj A equivalent}, Conjecture A would thus be implied by

\begin{eqn} \label{conj B'}
If $V < P$ is maximal, then every irreducible $\F_p[Out(P,V)]$--module occurs in $(\F_p[Aut(P,V)/Aut_0(P,V)] \otimes P_VH^*(C_P(V)))^{W_P(V)}$.
\end{eqn}

The $\F_p[Out(P,V)]$--module $(\F_p[Aut(P,V)/Aut_0(P,V)] \otimes P_VH^*(C_P(V)))^{W_P(V)}$ contains $\F_p[Aut(P,V)/Aut_0(P,V)]^{W_P(V)} \otimes P_VH^*(C_P(V))^{W_P(V)}$ as a submodule, and this submodule rewrites as
$$ \Ind_{Out_0(P,V)}^{Out(P,V)}(P_VH^*(C_P(V))^{W_P(V)}).$$

Conjecture B is the statement that every irreducible $\F_p[Out_0(P,V)]$--module occurs in $P_VH^*(C_P(V))^{W_P(V)}$. Thus Conjecture B would imply \ref{conj B'}, and we have proved

\begin{lem} \label{B implies A lemma}  If a pair $(P,V)$ satisfies Conjecture B, then $P$ satisfies Conjecture A.
\end{lem}

Now we prove that Conjecture C implies Conjecture B, establishing \thmref{conj C implies conj B thm}.
More precisely, let $V < P$ be maximal, and assume that every irreducible \\
$\F_2[Out_0(C_P(V),V)]$--module occurs as a submodule in $P_VH^*(C_P(V))$.  We will show that then every irreducible $\F_2[Out_0(P,V)]$--module occurs as a composition factor in $P_VH^*(C_P(V))^{W_P(V)}$.

The key here is again to use \corref{CP(V) prop}, our application of the Thompson $A \times B$ lemma, which tells us that the kernel of $Aut_0(P,V) \ra Aut_0(C_P(V),V)$ is a $2$--group.  Thus our assumption lets us conclude that every irreducible $\F_2[Aut_0(P,V)]$--module occurs as a submodule in $P_VH^*(C_P(V))$.

Now let $S$ be an irreducible $\F_2[Out_0(P,V)]$--module.  Equivalently, $S$ can be regarded as an irreducible $\F_2[Aut_0(P,V)]$--module fixed by the subgroup $W_P(V)$.  From our remarks in the last paragraph, there exists an inclusion $S \subset P_VH^*(C_P(V))$ of $\F_2[Aut_0(P,V)]$--modules, and thus an inclusion $$S = S^{W_P(V)} \subset P_VH^*(C_P(V))^{W_P(V)}$$ of $\F_2[Out_0(P,V)]$--modules.

\section{Some 2--central groups satisfying Conjecture C} \label{Conj C examples section}

If $P$ is a 2--central group, then $C= \Omega_1(P) < P$ is the unique maximal elementary subgroup and is central.  In this case, we write $Out_0(P)$ for $Out_0(P,C)$.

Recall that $P$ satisfies Conjecture C if every irreducible $\F_2[Out_0(P)]$--module occurs as a submodule of the algebra of primitive $P_CH^*(P)$.  In this section we verify this conjecture for many 2--central groups.

Our strategy is the following.  First, we verify the conjecture for various indecomposable 2--central groups, including all of order up to 64. Second, we find various side hypotheses on a pair $(P,Q)$ ensuring that if $P$ and $Q$ satisfy Conjecture C, so does $P \times Q$.

  All calculations have been done by hand, using information listed on the website \cite{carlson website} or the book \cite{carlson et al}.  For starters, we note that a 2--group $P$ is 2--central exactly when $\rk (Z(P)) = \rk(P)$, so 2--central groups can be easily identified. We occasionally use notation like `32\#10' to denote a 2--group: this would be the group of order 32 and numbered 10  on \cite{carlson website} or in  \cite{carlson et al}.

\subsection{Indecomposable 2--central groups of order dividing 64}

Let $P$ be 2--central. Recall that $P_CH^0(P) = \F_2$, the trivial module, so $P$ rather trivially satisfies Conjecture C whenever $Out_0(P)$ is a 2--group.  This is often the case.

\begin{lem}  Let $P$ be an indecomposable 2--central groups $P$ of order dividing 64.  Then  $Out_0(P)$ is a 2--group unless $P$ is one of the three groups 8\#5 ($Q_8$), 64\#162, or 64\#187 ($G_2$ as in \exref{2C examples}).
\end{lem}

We sketch how one checks this.  Firstly, here is a table of the indecomposable 2--central groups $P$ having {\em any} nontrivial automorphisms of order prime to 2, and the order of $Out(P)$ modulo its 2--Sylow subgroup.

\bigskip

\begin{center}
\begin{tabular}{|c|c|} \hline
$P$ & $2^{\prime}$ part of $|Out(P)|$  \\   \hline
8\#5  & 3 \\ \hline
32\#18 & 3 \\ \hline
64\#30 & 3 \\ \hline
64\#82 & 3 \\ \hline
64\#93 & 3 \\ \hline
64\#145 & 3 \\ \hline
64\#153 & 21 \\ \hline
64\#162 & 3 \\ \hline
64\#187 & 15 \\ \hline
\end{tabular}
\end{center}

\bigskip

\noindent We note that 8\#5 is $Q_8$, 64\#153 is the 2--Sylow subgroup of the Suzuki group $Sz(8)$, and 64\#187 is the 2--Sylow subgroup of $SU_3(\F_{4})$ as in \exref{2C examples}.

Now one checks that, except when $P$ is 8\#5 or 64\#162, an automorphism of $P$ of order 3 does {\em not} fix $\Omega_1(P)$, and similarly an automorphism of 64\#153 of order 7 does not fix $\Omega_1(P)$.

\begin{lem} The three groups 8\#5 ($Q_8$), 64\#162, and 64\#187 ($G_2$) all satisfy Conjecture C.
\end{lem}

To see this, we look at each group in turn.

If $P = Q_8$ with center $C$, then $Out_0(P)$ is the symmetic group of order 6, and there is a 2 dimensional $\F_2[Out_0(P)]$ irreducible $S$.  As $\F_2[Out_0(P)]$--modules, $S \simeq H^1(P/\Phi(P)) = H^1(P/C) = P_CH^1(P)$, and we are done.

If $P$ is 64\#162 with center $C$, the order of $Out_0(P)$ has the form $3\cdot 2^a$.  It follows that there must be exactly one irreducible $\F_2[Out_0(P)]$--module $S$ on which an element of order 3 acts nontrivially.  $H^1(P/\Phi(P))=H^1(P/C) = P_CH^1(P) = H^1(P)$ is 4 dimensional, with basis $\{z,y,x,w\}$ as in \cite{carlson et al}.  Examining the action of the generating automorphisms on $H^1(P)$, one sees that the 2 dimensional subspace with basis $\{x+z,w\}$ is an $Out_0(P)$--submodule acted on nontrivially by an automorphism of order 3.  Thus this submodule must be isomorphic to $S$, and we are done.

If $P$ is 64\#187 with center $C$, the order of $Out_0(P)$ has the form $5\cdot 2^a$.  It follows that there must be exactly one irreducible $\F_2[Out_0(P)]$--module $S$ on which an element of order 5 acts nontrivially. $H^1(P/\Phi(P))=H^1(P/C) = P_CH^1(P) = H^1(P)$ is four dimensional, and is an irreducible $\F_2[\Z/5]$--module as seen in \secref{2C example section}.  Thus $P_CH^1(P) \simeq S$, and we are done.

In summary, in this subsection we have checked

\begin{prop} \label{indec conj C prop} If $P$ is an indecomposable 2--central subgroup of order dividing 64, with $C=\Omega_1(P)$, then $|Out_0(P)|$ has the form $2^a$, $3\cdot 2^a$, or $5 \cdot 2^a$, and all irreducible $\F_2[Out_0(P)]$--modules occur as submodules of $P_CH^0(P)=\F_2$ or $P_CH^1(P)= H^1(P/C)$.
\end{prop}

\subsection{On $Out_0(P\times Q)$}

The following Kunneth formula for primitives is easily verified.

\begin{lem}  If $C < P$ and $D < Q$ are central elementary abelian subgroups, then
$$ P_{C\times D}H^*(P\times Q) = P_CH^*(P) \otimes P_DH^*(Q).$$
\end{lem}

\begin{defn}  Call a pair $(P,Q)$ of 2--central groups {\em good} if the algebra homomorphism
$$ \F_2[Out_0(P) \times Out_0(Q)] \ra \F_2[Out_0(P \times Q)]$$
induces an isomorphism on their semisimple quotients.
\end{defn}

In general, given two finite groups $G$ and $H$, the irreducible $\F_2[G \times H]$ modules will be direct summands of modules of the form $S \otimes T$, where $S$ is an irreducible $\F_2[G]$--module, and $T$ is an irreducible $\F_2[H]$--module.  Thus the lemma implies

\begin{cor}  If a pair $(P,Q)$ of 2--central groups is good, and both $P$ and $Q$ satisfy Conjecture C, so does $P \times Q$.
\end{cor}

Now we aim to prove that pairs of 2--central groups are frequently good.  The following is our key observation.

\begin{lem}  Let $P$ and $Q$ be 2--central.  Given
$$ \alpha = \left(
\begin{array}{cc} \alpha_P & \alpha_{P,Q} \\ \alpha_{Q,P} & \alpha_Q \end{array}\right): P \times Q \ra P \times Q$$
in $Aut_0(P \times Q )$, the following conditions on the component homomorphisms must hold. \\

\noindent{\bf (a)} $\Omega_1(P) \subseteq Ker( \alpha_{Q,P})$, and  $\Omega_1(Q) \subseteq Ker( \alpha_{P,Q})$. \\

\noindent{\bf (b)} $Im (\alpha_{Q,P}) \subseteq Z(Q)$, and $Im (\alpha_{P,Q}) \subseteq Z(P)$. \\

\noindent{\bf (c)} $\alpha_P \in Aut_0(P)$, and $\alpha_Q \in Aut_0(Q)$.
\end{lem}

\begin{proof}  As $\alpha$ must be the identity when restricted to $\Omega_1(P \times Q) = \Omega_1(P) \times \Omega_1(Q)$, condition (a) is clear. Furthermore, $\alpha_P$ and $\alpha_Q$ must be the identity when respectively restricted to $\Omega_1(P)$ and $\Omega_1(Q)$.

As $\alpha$ is a homomorphism, it follows that $Im (\alpha_{Q,P})$ and $Im(\alpha_Q)$ must commute, and similarly $Im (\alpha_{P,Q})$ and $Im(\alpha_P)$ must commute.

It follows that both (b) and (c) will follow if we can show that $\alpha_P$ and $\alpha_Q$ are isomorphisms.

We are assuming that $\alpha$ is invertible.  Thus there exists
$$ \beta = \left(
\begin{array}{cc} \beta_P & \beta_{P,Q} \\ \beta_{Q,P} & \beta_Q \end{array}\right): P \times Q \ra P \times Q$$
such that $\beta \circ \alpha: P \times Q \ra P \times Q$ is the identity.  In particular,
\begin{eqn} \label{id condition} The composite $P \xra{\alpha_P, \alpha_{Q,P}} P \times Q \xra{\beta_P \cdot \beta_{P,Q}} P$ is the identity.
\end{eqn}

We will show that $\beta_P \circ \alpha_P: P \ra P$ is epic.  Let $\gamma = \beta_{P,Q} \circ \alpha_{Q,P}: P \ra P$.  Then \ref{id condition} says
$$ x = \beta_P(\alpha_P(x))\gamma(x) \text{ for all } x \in P.$$
By condition (a), $\Omega_1(P) \subseteq Ker( \gamma)$, and so the order of $\gamma(x)$ will be strictly less than the order of $x$ for all $x \neq e \in P$.  Thus we can prove that $\beta_P \circ \alpha_P$ is epic by induction on the order of $x$, as follows. Given $x \in P$, there exists $z$ such that $\beta_P(\alpha_P(z)) = \gamma(x)$.  But then
$$\beta_P(\alpha_P(xz)) = \beta_P(\alpha_P(x))\beta_P(\alpha_P(z)) = \beta_P(\alpha_P(x))\gamma(x) = x.$$
\end{proof}

\begin{thm} A pair $(P,Q)$ of 2--central groups is good if every homomorphism
$P/\Omega_1(P) \ra Z(Q)$ has image contained in $\Phi(Q)$.
\end{thm}

\begin{proof}  Let $\widetilde{Aut}_0(P)$ denote the image of $Aut_0(P) \ra Aut(P/\Phi(P))$.  \propref{phi(P) prop} implies that $\F_2[Aut_0(P)] \ra \F_2[\widetilde{Aut}_0(P)]$ induces an isomorphism on semisimple quotients, and similarly with $P$ replaced by $Q$ and $P \times Q$.

Using parts (a) and (b) of the lemma, the hypothesis that every homomorphism $P/\Omega_1(P) \ra Z(Q)$ has image contained in $\Phi(Q)$ implies that for all $\alpha \in Aut_0(P\times Q)$, the component $\alpha_{Q,P}: P \ra Q$ induces the zero map $P/\Phi(P) \ra Q/\Phi(Q)$.  Again using the lemma, it follows that
$$ \widetilde{Aut}_0(P \times Q) = \left(
\begin{array}{cc} \widetilde{Aut}_0(P) & B \\ 0 & \widetilde{Aut}_0(Q) \end{array}
\right),$$
where $B$ is a quotient of $\Hom(Q/\Omega_1(Q), Z(P))$, and is thus a 2--group.

The inclusion $\widetilde{Aut}_0(P) \times \widetilde{Aut}_0(Q) \ra \widetilde{Aut}_0(P\times Q)$ thus has a retract with kernel a 2--group, and so induces an isomorphism on the semisimple quotients of the associated $\F_2$ group rings.
\end{proof}

\begin{cor} \label{good pair cor} A pair $(P,Q)$ of 2-central groups is good if any of the following conditions hold. \\

\noindent{\bf (a)} $Z(Q) \subseteq \Phi(Q)$. \\

\noindent{\bf (b)} $P$ is elementary abelian. \\

\noindent{\bf (c)} $P/[P,P]\Omega_1(P)$ is elementary abelian, and $Q$ has no $\Z/2$ summands. \\

\noindent{\bf (d)} $Q = \Z/2^t$ and $P$ has exponent dividing $2^t$.
\end{cor}

In summary, we have shown

\begin{prop} If $P$ and $Q$ satisfy any of the conditions in the last corollary, and both satisfy Conjecture C, then $P \times Q$ also satisfies Conjecture C.
\end{prop}

\subsection{Many 2--central groups satisfy Conjecture C}

When one does a census of indecomposable noncyclic 2-central groups, one finds 1 of order 8, 2 of order 16, 9 of order 32, and 41 of order 64.\footnote{Group 32\#25 is mistakenly listed at indecomposable in both \cite{carlson website} and \cite{carlson et al}.  It is isomorphic to $\Z/2 \times Q_{16}$.}  Condition (a) of \corref{good pair cor} holds for all but one of these groups.

\begin{prop} Let $Q$ be an indecomposable 2--central group of order dividing 64.  Then $Z(Q) \subseteq \Phi(Q)$ unless $Q$ is cyclic or 64\#30.  The group 64\#30 has exponent 4, and thus satisfies the condition for $P$ listed in (c) of \corref{good pair cor}.
\end{prop}

Checking this is made easier by the observation that a 2-central group $Q$ has no $\Z/2$ summands exactly when $\Omega_1(Q) \subseteq \Phi(Q)$.  Thus $Z(Q) \subseteq \Phi(Q)$ whenever $Z(Q)$ is elementary abelian. This is the case for all but 11 of the 53 indecomposable noncyclic 2-central groups in question.

We can now easily prove \thmref{2 central conj thm}, which we restate here.

\begin{thm} \ Conjecture C is true for all 2--central groups that can be written as the product of groups of order dividing 64.
\end{thm}

\begin{proof}  By \propref{indec conj C prop}, all indecomposable 2--central groups of order dividing 64 satisfy Conjecture C.  Now suppose a 2--central group $P$ is a product $P = A \times B \times C$ where $A$ is abelian with no $\Z/2$ summands, $B$ is a product of noncyclic indecomposable summands of order dividing 64, and $C$ is elementary abelian.  Repeatedly using condition (d) of \corref{good pair cor}, one deduces that $A$ satisfies Conjecture C.\footnote{Indeed $Aut_0(A)$ is a 2--group, as is well known.}  Repeatedly using conditions (a) and (c), one then deduces that $A \times B$ satisfies Conjecture C, and finally that $P = A \times B \times C$ does, thanks to condition (b).
\end{proof}

\section{A discussion of \exref{last group example}} \label{examples section}

In this section we discuss in detail \exref{last group example}, and use this to illustrate \thmref{R_d theorem} and its proof.  So let $P$ be the group 64\#108: all of our information about this group is available on \cite{carlson website}.

\subsection{Subgroup structure}  The commutator subgroup $Z=[P,P]$ has order 2. The center $C$ is elementary abelian of rank 2, and $C = \Phi(P)$, so $Z < C$ and $P/C$ is elementary abelian of rank 4.   There is a unique maximal elementary abelian group $V$ of rank 3, and its centralizer $K=C_P(V)$ has order 32, so that $W_P(V) = G/K \simeq \Z/2$.  More precisely, the 2--central group $K$ is isomorphic to $(\Z/2)^2  \times Q_8$, with $Q_8$ embedded so that $V \cap Q_8 = Z$.\footnote{In terms of the generators $g_1, \dots, g_6$ on the website, $Z=\langle g_6 \rangle, C=\langle g_3,g_6 \rangle$, $V=\langle g_2,g_3,g_6 \rangle$, $Q_8=\langle g_4,g_5,g_6 \rangle$, and $K = \langle g_2,g_3,g_4,g_5,g_6 \rangle$.  $K$ is maximal subgroup \#11.}

We have the following picture of $\A^C(P)$:

\begin{equation}
\xymatrix{
C  \ar[r] & V \ar@(ur,dr)^{\Z/2}   }
\end{equation}
and from this it is already clear that $\HH^*(P) = H^*(V)^{\Z/2}$.

\subsection{$H_1(P)$ as an $\F_2[Out(P)]$--module}  The group $Out(P)$ has order $768 = 3 \cdot 2^8$.  From this, one can deduce that $\F_2[Out(P)]$ has precisely two irreducible modules, the trivial module `1', and another one which we will call `$S$' of dimension at least two.  Furthermore, $S$ occurs as a composition factor in an $\F_2[Out(P)]$--module $M$ if and only if any element of order 3 in $Out(P)$ acts nontrivially on $M$.

Now we consider $H_1(P)=P/C$ as an $\F_2[Out(P)]$--module.  It is nontrivial by \propref{phi(P) prop}, thus S occurs as a composition factor.  Notice that each of the subgroups $C$, $V$, and $K$ are characteristic, and recall that $Inn(P)$ acts trivially on $H_1(P)$.  Thus there is a composition series of $\F_2[Out(P)]$--modules:
$$   V/C \subset K/C \subset P/C,$$
with composition factors $V/C$, $K/V$, and $P/K$.  By dimension considerations, $V/C$ and $P/K$ are copies of the trivial module, and thus $K/V$ must be $S$, which is thus two dimensional.

\subsection{$H^*(V)$ as an $\F_2[Aut(P)]$--module}  An automorphism of $P$ order 3 fixes $V$.\footnote{See automorphism \#7 on the Carlson website.}  Thus, as an $Aut(P)$--module, $H^*(V)$ has only trivial composition factors.

\subsection{ $P_VH^*(K)$ and $P_CH^*(K)$}  We have maps of unstable algebras equipped with with $Aut(P)$ action:
\begin{equation} P_VH^*(K)\hookrightarrow P_CH^*(K) \xla{j^*} P_CH^*(P),
\end{equation}
where $j:K \ra P$ is the inclusion. In degree 1, this reads
$$ \Hom(K/V, \F_2) \hookrightarrow \Hom(K/C, \F_2) \la \Hom(P/C, \F_2),$$
and so we see that $j^*$ is onto in degree 1, and that $P_VH^1(K)$ is a copy of $S$.

The maps of pairs $(Q_8,Z) \ra ((\Z/2)^2 \times Q_8, (\Z/2)^2 \times Z) = (K,V)$ induces an isomorphism of algebras:
$$P_VH^*(K) \simeq P_ZH^*(Q_8).$$

The algebra $P_ZH^*(Q_8)$ is familiar: the calculation of $H^*(Q_8)$ using the Serre spectral sequence associated to $Z \ra Q_8 \ra Q_8/Z$ reveals that $P_ZH^*(Q_8) = \text{Im }\{ H^*(Q_8/Z) \ra H^*(Q_8)\} = B^*$, where $B^*$ is the Poincar\'e duality algebra $\F_2[x,w]/(x^2 + xw + w^2, x^2w + xw^2)$, where $x$ and $w$ both have degree 1.  As an $\F_2[Aut(P)]$--module, $P_VH^*(K) = B^*$ is given by
\begin{equation*}
B^d =
\begin{cases}
1 & \text{if } d = 3 \\ S & \text{if } d = 2 \\ S & \text{if } d = 1 \\ 1 & \text{if } d = 0 .
\end{cases}
\end{equation*}

From this we learn that $P_CH^*(K) \simeq B^*[y]$ where $y$ is also in degree 1, and thus is generated by elements in degree 1.  It follows that $j^*: P_CH^*(P) \ra P_CH^*(K)$ is onto, and then that $Inn(P)$ acts trivially on both $P_VH^*(K)$ and $P_CH^*(K)$.

\subsection{$\bar R_dH^*(P)$}  Applying the formula for $\bar R_dH^*(P)$ in \propref{good Rd formula prop} to our group, we see that, for all $d$, there is a pullback diagram of graded $\F_2[Out(P)]$--modules:
\begin{equation*}
\xymatrix{
\bar R_dH^*(P) \ar[d] \ar[r] & H^*(C) \otimes P_CH^d(P) \ar[d]^{1\otimes j^*}  \\
(H^*(V) \otimes P_VH^d(K))^{Inn(P)} \ar[r] & H^*(C) \otimes P_CH^d(K). }
\end{equation*}

By our comments above, this simplifies: for all $d$, there is a pullback diagram of unstable $\F_2[Out(P)]$--modules:
\begin{equation*}
\xymatrix{
\bar R_dH^*(P) \ar[d] \ar[r] & H^*(C) \otimes P_CH^*(P) \ar[d]^{1\otimes j^*}  \\
H^*(V)^{\Z/2} \otimes P_VH^*(K) \ar[r] & H^*(C) \otimes P_CH^*(K). }
\end{equation*}

As $j^*$ is onto, we conclude
\begin{prop} \label{Rd calc} For all $d$, there is a short exact sequence of unstable $\F_2[Out(P)]$--modules $ 0 \ra H^*(C) \otimes \ker{j^d} \ra \bar R_dH^*(P) \ra H^*(V)^{\Z/2} \otimes P_VH^d(K) \ra 0,$
where $H^*(C)$ and $H^*(V)^{\Z/2}$ have only trivial composition factors, while $P_VH^*(K) \simeq B^*$.
\end{prop}

We note that in the notation of \secref{proof of the main theorem}, $\Ess^*(C) = \ker{j^*}$ and $\Ess^*(V) = P_VH^*(K)$.

From \propref{Rd calc}, one can already place the 3's occurring in the table in \exref{last group example}. The 2's will be determined by knowing the composition factors of the ideal $\ker j^* \subset P_CH^*(P)$.  This we discuss next.

\subsection{$P_CH^*(P)$ and the ideal $\ker j^*$}  Obviously, $\ker j^0 = 0$, and we also know that $\ker j^1$ is one dimensional.   Using more detailed information about $H^*(P)$, we outline how all of the composition factors of $\ker j^*$ can be determined.

$H^*(P)$ is generated by classes $z,y,x,w,v,u,t$ of degrees 1,1,1,1,2,5,8 respectively, and the kernel of the restriction $H^*(P) \ra H^*(K)$ is identified as the ideal $(z)$.\footnote{See the description of the depth essential cohomology on the website.}  Furthermore, this ideal is a free module on $\F_2[v,t]$ on an explicit basis of 24 monomials in the generators of degree 1.

All the degree 1 generators are $H^*(C)$--primitive, while $v$ and $t$ restrict to algebraically independent elements of $H^*(C)$.  It follows that $\ker{j^*} = P_CH^*(P) \cap (z)$ is precisely the $\F_2$--span of this basis of mononomials.

This allows us to conclude the following.

Firstly $P_CH^*(P) \subset H^*(P)$ is precisely the subalgebra generated by the classes in degree 1.  Explicitly, $P_CH^*(P)$ is
$$ A^* = \F_2[z,y,x,w]/(z^2,zy+x^2+xw+w^2,zy^2+x^2w+xw^2,y^2w^3+yw^4+w^5).$$
For our purposes, this can be simplified as follows.  From the second and third relations, one can deduce that $x^3 = z(yx+y^2)$ and $w^3 = z(yw+y^2)$.  Using these and the first relation, one deduces that $y^2w^3+yw^4+w^5 = zy^4$.  So we have
$$ A^* = \F_2[z,y,x,w]/(z^2,zy^4, zy+x^2+xw+w^2,zy^2+x^2w+xw^2).$$

Secondly, $\ker{j^*} = \Sigma (A^*/Ann(z))$, as unstable $\F_2[Out(P)]$--modules.  Since $z^2 = zy^4 = 0$, $A^*/\text{Ann}(z)$ is a quotient of $A^*/(z,y^4) = \F_2[y,x,w]/(y^4, x^2+xw+w^2,x^2w+xw^2) = B^*[y]/(y^4)$.
Inspection of the degrees of the monomial basis shows that $A^*/Ann(z)$ has dimension 1,3,5,6,5,3,1 in degrees 0,1,2,3,4,5,6.  This agrees with $B^*[y]/(y^4)$, and we can conclude that $Ann(z) = (z,y^4)$, and then

\begin{prop} $\ker{j^*} \simeq \Sigma (B^*[y]/(y^4))$ as unstable $\F_2[Out(P)]$--modules.
\end{prop}


\begin{thebibliography}{99}


% \bibitem[AKM]{akm} A.Adem, D.B.Karagueuzian, and J.Min\'a\v{c}, {\em On the cohomology of Galois groups
% determined by Witt rings},  Adv. Math.  {\bf 148}  (1999),  105--160.

% \bibitem[AP]{ap} A.Adem and J.Pakianathan, {\em On the cohomology of central Frattini extensions},  J.
% Pure Appl. Algebra  {\bf 159}  (2001), 1--14.


\bibitem[Alp]{alperin}  J.~ L.~ Alperin, {\em Local Representation Theory}, Cambridge Studies in Advanced Math. {\bf  11}, Cambridge U. Press, 1986.

\bibitem[Asch]{asch}  M.~ Aschbacher, {\em Finite Group Theory}, 2nd edition, Cambridge Studies in Advanced Math. {\bf  10}, Cambridge U. Press, 2000.


%% \bibitem[BC]{bc} D.J.Benson, J.F.Carlson, {\em Projective resolutions and Poincaré duality complexes},
%% Trans. Amer. Math. Soc. {\bf 342}  (1994),  447--488.

%% \bibitem[Be]{b} D.Benson, {\em Dickson invariants, regularity and computation in group cohomology},
%% Illinois J. Math.  {\bf 48}  (2004),  171--197.

\bibitem[BoZ]{bour zarati} D. Bourguiba and S. Zarati, {\em Depth and the Steenrod algebra}. With an appendix by J. Lannes.  Invent. Math. {\bf 128} (1997), 589--602.

\bibitem[BrH]{broto henn} C.Broto and H.-W. Henn, {\em Some remarks on central elementary abelian $p$--subgroups and cohomology of classifying spaces},  Quart. J. Math. {\bf 44}  (1993),  155--163.

\bibitem[BrZ1]{broto zarati 1} C.Broto and S.Zarati, {\em Nil--localization of unstable algebras over the Steenrod algebra},  Math. Zeit. {\bf 199}  (1988),  525--537.

\bibitem[BrZ2]{broto zarati 2} C.Broto and S.Zarati, {\em On sub--$\A_p^*$--algebras of $H^*(V)$},  Springer L. N. Math.  {\bf 1509}  (1992),  35--49.

%% \bibitem[BP]{bp} W.Browder and J. Pakianathan, {\em Cohomology of uniformly powerful $p$-groups},
%% Trans. Amer. Math. Soc.  {\bf 352}  (2000),  2659--2688.

\bibitem[Ca]{carlson website} Jon Carlson, {\em Mod 2 cohomology of 2 groups}, MAGMA computer computations on the website http://www.math.uga.edu/{$\sim$}lvalero/cohointro.html.

\bibitem[CTVZ]{carlson et al} J.F.Carlson, L.Townsley, L.Valeri-Elizondo,and M. Zhang, {\em Cohomology rings of finite groups}. With an appendix, {\em Calculations of cohomology rings of groups of order dividing 64}, by Carlson, Valeri-Elizondo and Zhang. Algebras and Applications {\bf 3}, Kluwer, Dordrecht, 2003.

\bibitem[Cr]{crabb} M.C.Crabb, {\em Dickson-Mui invariants}, Bull. London Math. Soc. {\bf 37} (2005), 846--856.

\bibitem[D]{dickson} L.E.Dickson, {\em A fundamental system of invariants of the general modular linear group with a solution of the form problem},  T..A.M.S. {\bf 12}  (1911),  75--98.

\bibitem[DS]{ds} T. Diethelm and U. Stammbach,  {\em On the module structure of the mod {$p$} cohomology of a {$p$}-group},  Arch. Math.  {\bf 43}  (1984),  488--492.

\bibitem[D]{duflot} J. Duflot,  {\em Depth and equivariant cohomology}, Comm. Math. Helv. {\bf 56}  (1981),  627--637.

%% \bibitem[F-A]{Fernandez}  G.~A.~ Fern\'andez-Alcober, {\em An introduction to finite $p$-groups: regular %% $p$-groups and groups of maximal class},  Matem\'atica Contempor\^anea {\bf 20} (2001), 155--226.


\bibitem[Gor]{gor}  D. Gorenstein, {\em Finite Groups}, 2nd edition, Chelsea Publishing, New York, 1980.

\bibitem[Gr]{green}  D.J.Green, {\em The essential ideal in group cohomology does not square to zero},  J. Pure Appl. Algebra {\bf 193}  (2004),  no. 1-3, 129--139.


\bibitem[HK]{hk} J.C.Harris and N.J.Kuhn, {\em Stable decompositions of classifying spaces of finite abelian p-groups}, Math. Proc. Camb. Phil. Soc. \  {\bf 103}(1988), 427--449.

\bibitem[H]{henn}    H.-W. Henn, {\em Finiteness properties of injective resolutions of certain unstable modules over the Steenrod algebra and applications}, Math. Ann. {\bf 291} (1991), 191-203.

\bibitem[HLS1]{hls1}  H.-W. Henn, J. Lannes, and L. Schwartz, {\em The categories of unstable modules and  unstable algebras modulo nilpotent objects}, Amer. J. Math. {\bf 115} (1993), 1053-1106.

\bibitem[HLS2]{hls2} H.--W.Henn, J.Lannes, and L.Schwartz, {\em Localizations of unstable $A$-modules and equivariant mod $p$ cohomology},  Math. Ann. {\bf 301}  (1995), 23--68.


\bibitem[HP]{henn priddy} H.--W.Henn and S.Priddy, {\em $p$--nilpotence, classifying space indecomposability, and other properties of almost all finite groups}, Comm. Math. Helv. \  {\bf 69}(1994), 335--350.
\bibitem[Hi]{higginbottom} R.~ Higginbottom, Ph.D. ~thesis, University of Virginia, 2005.

\bibitem[K1]{k1} N.J.Kuhn, {\em Character rings in algebraic topology}, Advances in Homotopy Theory (Cortona 1988), L. M. S. Lectures Notes \ {\bf 139}(1989), 111--126.


\bibitem[K2]{genrep1} N.J.Kuhn, {\em Generic representations of the finite general linear groups and the Steenrod algebra:I}, Amer. J. Math.  {\bf 116}(1994), 327--360.

\bibitem[K3]{k2} N.J.Kuhn, {\em On topologically realizing modules over the Steenrod algebra}, Ann. Math. {\bf 141}(1995), 321--347.

\bibitem[K4]{k3} N.J.Kuhn, {\em Cohomology primitives associated to central extensions}, Oberwolfach Reports {\bf 2} (2005), 2383--2386.

\bibitem[K5]{k4} N.J.Kuhn, {\em Cohomology primitives associated to central extensions}, in preparation.

\bibitem[LZ]{lannes zarati} J.Lannes, and S.Zarati, {\em Sur les $\U$--injectifs},  Ann. Sci. Ec. Norm. Sup. {\bf 19}  (1986), 303--333.

\bibitem[MP]{mp} J.Martino and S.Priddy, {\em On the dimension theory of dominant summands},  Adams Memorial Symposium on Algebraic Topology, 1 (Manchester, 1990),  L.M.S. Lect. Note Ser. {\bf 175} (1992), 281--292.

\bibitem[Mat]{matsumura} H.Matsumura, {Commutative algebra}, 2nd edition, Math. Lect. Note. Series, Benjamin, 1980.

%% \bibitem[MiS]{} J.Min\'a\v{c} and M.Spira, {\em Witt rings and Galois groups}, Ann. Math. {\bf 144}
%% (1996), 35--60.

\bibitem[N]{nishida} G. Nishida, {\em Stable homotopy type of classifying spaces of finite groups},  Algebraic and topological theories (Kinosaki 1984), Kinokuniya, Tokyo, 1986, 391--404.

\bibitem[Q1]{quillen} D. Quillen, {\em The spectrum of an equivariant cohomology ring {I}},
Ann. Math. {\bf 94} (1971), 549-572.

\bibitem[Q2]{quillen part II} D. Quillen, {\em The spectrum of an equivariant cohomology ring {II}},
Ann. Math. {\bf 94} (1971), 573-602.

%% \bibitem[Q3]{quillen2} D. Quillen, {\em The ${\rm mod}$ $2$ cohomology rings of extra-special $2$-groups %% and the {S}pinor groups}, Math. Ann. {\bf 194} (1971), 197-212.


\bibitem[S1]{s1}   L. Schwartz, {\em La filtration nilpotente de la cat\'egorie $\U$ et la cohomologie des espaces de lacets} , Algebraic Topology -- Rational Homotopy (Louvain la Neuve 1986), S. L. N. M. 1318 (1988), 208-218.

\bibitem[S2]{s2}  L. Schwartz, {\em Modules over the Steenrod algebra and Sullivan's fixed point conjecture}, Chicago Lectures in Math., U. Chicago Press, 1994.

\bibitem[Sy]{symonds} P. Symonds, {\em The action of automorphisms on the cohomology of a $p$-group}, Math. Proc. Camb. Phil. Soc. {\bf 127} (1999), 495--496.

\end{thebibliography}
\end{document}